\documentclass{amsart}
\usepackage{amsthm,amsmath,amsfonts}
\usepackage{mathrsfs}
\usepackage[colorlinks=true,linktocpage=true]{hyperref}
\usepackage{accents}
\usepackage{bbm,dsfont}
\usepackage{hyperref}
\hypersetup{urlcolor=blue, citecolor=blue, linkcolor=red}
\usepackage{framed}
\usepackage{graphics}
\usepackage{graphicx}
\usepackage{color}
\usepackage{tocvsec2}
\numberwithin{equation}{section}
\numberwithin{table}{section}
\font\tenscrpt=eusm10
\font\sevenscrpt=eusm10 scaled 700
\font\fivescrpt=eusm10 scaled 500
\newfam\eusmfam
\textfont\eusmfam=\tenscrpt
\scriptfont\eusmfam=\sevenscrpt
\scriptscriptfont\eusmfam=\fivescrpt

\newtheorem{thm}{Theorem}[section]
\newtheorem{cor}{Corollary}[section]
\newtheorem{lem}{Lemma}[section]
\newtheorem{prop}{Proposition}[section]

\theoremstyle{definition}
\newtheorem{defn}{Definition}[section]

\newtheorem{rem}{Remark}[section]
\newtheorem{notn}{Notation}[section]
\newcommand{\thmref}[1]{Theorem~\ref{#1}}
\newcommand{\secref}[1]{Section~\ref{#1}}

\newcommand{\lemref}[1]{Lemma~\ref{#1}}
\newcommand{\coref}[1]{Corollary~\ref{#1}}

\newcommand{\defnref}[1]{Definition~\ref{#1}}
\newcommand{\remref}[1]{Remark~\ref{#1}}
\newcommand{\eqnref}[1]{{\rm (\ref{#1})}}

\def\qed{\quad\vcenter{\hrule\hbox{\vrule height.6em\kern.6em\vrule}\hrule}}
\newenvironment{pf}{{\bigskip\textit{\newline Proof.}\quad}}{$\qed$\bigskip\newline}
\newenvironment{pf*}[1]{{\bigskip\textit{\newline#1.}\quad}}{$\qed$\bigskip\newline}
\def\ds{\displaystyle}
\def\ts{\textstyle}

\def\i{\mathbf i}

\def\ptzs{{K}^{\text{\tiny{\sc{BM}}}}_{t;0,s} }

\def\Mpisxy {{K}^{\text{\tiny{\sc{ASP}}}^d}_{\i s;x,y}}

\def\psxz{{K}^{\text{\tiny{\sc{BM}}}^d}_{s;x}}

\def\ptsz{{K}^{\text{\tiny{\sc{BM}}}}_{t;s}}
\def\ptzs{{K}^{\text{\tiny{\sc{BM}}}}_{t;0,s}}
\def\peptsz{{K}^{\text{\tiny{\sc{BM}}}}_{\vep t;s}}

\def\K{{\mathbb K}}
\def\KBtxy{{\K}^{\text{\tiny{\sc{BTBM}}}^d}_{t;x,y}}
\def\KKStxy{{\K}^{\text{\tiny{\sc{LKS}}}^d}_{t;x,y}}

\def\KKSeptsxy{{\K}^{\text{\tiny{\sc{LKS}}}^d}_{\vepo(t-s);x,y}}

\def\KKSsxpz{{\K}^{\text{\tiny{\sc{LKS}}}^d}_{s;x+z}}
\def\FKKSsxipz{\hat{\K}^{\text{\tiny{\sc{LKS}}}^d}_{s;\xi+z}}

\def\KKStsx{{\K}^{\text{\tiny{\sc{LKS}}}^d}_{t-s;x}}
\def\KKStx{{\K}^{\text{\tiny{\sc{LKS}}}^d}_{t;x}}
\def\KKSeptx{{\K}^{\text{\tiny{\sc{LKS}}}^d}_{\vepo t;x}}
\def\KKSepotx{{\K}^{(\vepo,1)\text{\tiny{\sc{LKS}}}^d}_{t;x}}
\def\KKSepthtx{{\K}^{(\vep,\vth)\text{\tiny{\sc{LKS}}}^d}_{t;x}}

\def\KKSepthootx{{\K}^{(1,1)\text{\tiny{\sc{LKS}}}^d}_{t;x}}
\def\KKSepthoztx{{\K}^{(1,0)\text{\tiny{\sc{LKS}}}^d}_{t;x}}
\def\KKSepthtxy{{\K}^{(\vep,\vth)\text{\tiny{\sc{LKS}}}^d}_{t;x,y}}

\def\Ksfotx{{\K}^{\text{\tiny{\sc{SFO}}}^d}_{t;x}}
\def\FKsfotxi{\hat{\K}^{\text{\tiny{\sc{SFO}}}^d}_{t;\xi}}
\def\FKsfosxi{\hat{\K}^{\text{\tiny{\sc{SFO}}}^d}_{s;\xi}}
\def\Ksfosx{{\K}^{\text{\tiny{\sc{SFO}}}^d}_{s;x}}

\def\FKKStxi{\hat{\K}^{\text{\tiny{\sc{LKS}}}^d}_{t;\xi}}
\def\FKKSepthtxi{\hat{\K}^{(\vep,\vth)\text{\tiny{\sc{LKS}}}^d}_{ t;\xi}}
\def\FKKSsxi{\hat{\K}^{\text{\tiny{\sc{LKS}}}^d}_{s;\xi}}
\def\KKSrsx{{\K}^{\text{\tiny{\sc{LKS}}}^d}_{r-s;x}}
\def\KKSsptmrx{{\K}^{{\text{\tiny{\sc{LKS}}}^d}}_{s+(t-r);x}}
\def\KKSsx{{\K}^{{\text{\tiny{\sc{LKS}}}^d}}_{s;x}}
\def\KKSepsx{{\K}^{{\text{\tiny{\sc{LKS}}}^d}}_{\vepo s;x}}
\def\FKKSsptmrxi{\hat{\K}^{{\text{\tiny{\sc{LKS}}}^d}}_{s+(t-r);\xi}}
\def\FKKSsxi{\hat{\K}^{{\text{\tiny{\sc{LKS}}}^d}}_{s;\xi}}

\def\KBtxy{{\K}^{\text{\tiny{\sc{BTBM}}}^d}_{t;x,y}}
\def\KBtx{{\K}^{\text{\tiny{\sc{BTBM}}}^d}_{t;x}}
\def\FKBtxi{\hat{\K}^{\text{\tiny{\sc{BTBM}}}^d}_{t;\xi}}
\def\KKStxy{{\K}^{\text{\tiny{\sc{LKS}}}^d}_{t;x,y}}

\def\KKSepthtsxy{{\K}^{(\vep,\vth)\text{\tiny{\sc{LKS}}}^d}_{t-s;x,y}}

\def\utx{U(t,x)}

\def\usy{U(s,y)}

\def\unx{u_0(x)}

\def\uny{u_0(y)}
\def\un{u_0}

\def\sW{\mathscr W}
\def\sWi{{\mathscr W}^{(i)}}

\def\cC{\mathbb R}
\def\dD{\mathbb D}
\def\P{{\mathbb P}}

\def\EP{{\mathbb E}_{\P}}

\def\E{{\mathbb E}}
\def\N{{\mathbb N}}

\def\Rd{{\mathbb R}^{d}}

\def\R{\mathbb R}
\def\Rs{{\mathbb R}^2}
\def\B{{\mathbb B}}
\def\S{\mathbb S}
\def\T{\mathbb T}

\def\Rp{{\R}_+}

\def\Rd{{\mathbb R}^{d}}

\def\sF{{\mathscr F}}

\def\sFt{{\mathscr F}_t}

\def\OFFtP{(\Omega,\sF,\{\sFt\},\P)}
\def\OFFtPi{(\Omega^{i},\sF^{i},\{{\sFt}^{i}\},\P^{i})}

\def\C{\mathrm C}
\def\H{\mathrm H}

\def\eqdef{:=}

\def\D{\Delta}
\def\lap{\Delta}

\def\df#1#2{\ds{\frac{#1}{#2}}}
\def\tf#1#2{\ts{\frac{#1}{#2}}}
\def\lbl#1{\label{#1}}

\def\intrd{\int_{\Rd}}

\def\intrs{\int_{\Rs}}

\def\intrdzt{\int_{\Rd}\int_0^t}

\def\pa{\partial}

\def\lab{\left|}
\def\rab{\right|}
\def\lpa{\left(}
\def\rpa{\right)}
\def\lbk{\left[}
\def\rbk{\right]}
\def\lbr{\left\{}
\def\rbr{\right\}}
\def\bdf{\begin{defn}}
\def\edf{\end{defn}}
\def\bcr{\begin{cor}}
\def\ecr{\end{cor}}
\def\bnt{\begin{notn}}
\def\ent{\end{notn}}
\def\brm{\begin{rem}}
\def\erm{\end{rem}}
\def\blm{\begin{lem}}
\def\elm{\end{lem}}
\def\bpf{\begin{pf}}
\def\bpfs{\begin{pf*}}
\def\epf{\end{pf}}
\def\epfs{\end{pf*}}

\def\beq{\begin{equation}}
\def\beqs{\begin{equation*}}
\def\eeq{\end{equation}}
\def\eeqs{\end{equation*}}
\def\bsp{\begin{split}}
\def\esp{\end{split}}
\def\bc{\begin{cases}}
\def\ec{\end{cases}}
\def\bt{\begin{tabular}}
\def\et{\end{tabular}}
\def\ben{\begin{enumerate}}

\def\rencomrom{\renewcommand{\labelenumi}{(\roman{enumi})}}
\def\renrom{\renewcommand{\labelenumi}{(\roman{enumi})}}
\def\een{\end{enumerate}}

\def\bthm{\begin{thm}}
\def\ethm{\end{thm}}
\def\bpr{\begin{prop}}
\def\epr{\end{prop}}
\def\bcor{\begin{cor}}
\def\ecor{\end{cor}}
\def\bfr{\begin{framed}}
\def\efr{\end{framed}}
\def\ig{\iffalse}
\def\bcm{\iffalse}
\def\babs{\begin{abstract}}
\def\eabs{\end{abstract}}

\def\e{\mathrm{e}}

\def\vep{\varepsilon}
\def\vepo{\vep_{1}}
\def\vept{\vep_{2}}
\def\vth{\vartheta}
\def\vph{\varphi}
\def\pa{\partial}
\def\pat{\pa_{t}}
\def\eLKSabun{\e^{(\vep,\vth)}_{\mbox{\tiny LKS}}(a,b,\un)}
\def\eLKSazun{\e^{(\vep,\vth)}_{\mbox{\tiny LKS}}(a,0,\un)}
\def\sB{\mathscr{B}}
\def\Oo{\Omega^{(1)}}
\def\Ot{\Omega^{(2)}}
\def\Pii{\P^{(i)}}
\def\Ui{U^{(i)}}
\def\Uo{U^{(1)}}
\def\Ut{U^{(2)}}
\def\LawUP{\mathscr{L}_{\P}^{U}}
\def\LawVQ{\mathscr{L}_{\Q}^{U}}
\def\LawUiPi{\mathscr{L}_{\Pii}^{\Ui}}
\def\RU{R_{U}}
\def\RV{R_{V}}
\def\Q{\mathbb Q}

\def\Pii{\P^{(i)}}
\def\Pni{\P_n^{(i)}}
\def\Ptn{\tilde{\P}_n}
\def\Piti{\tilde{\P}^{(i)}}
\def\Poti{\tilde{\P}^{(1)}}
\def\Ptti{\tilde{\P}^{(2)}}
\def\Pno{\P_n^{(1)}}
\def\Pnt{\P_n^{(2)}}
\def\EP{{\Bbb E}_{\P}}
\def\EPno{{\Bbb E}_{\Pno}}
\def\EPnt{{\Bbb E}_{\Pnt}}
\def\Ft{{\mathscr{F}}_t}

\def\FT{{\mathscr{F}}_T}
\def\FTi{{\mathscr{F}}_T^{(i)}}
\def\filpspace{(\Omega, \mathscr{F}, \{{\mathscr{F}}_t\},\P)}

\def\filpspaceTni{(\Omega^{(i)}, \mathscr{F}_T^{(i)}, \{{\mathscr{F}}_t^{(i)}\},{\P}_n^{(i)})}

\def\filpspaceiti{(\Omega^{(i)}, {\mathscr{F}}^{(i)}, \{{\mathscr{F}}_t^{(i)}\},\Piti)}
\def\filpspace{(\Ot, \mathscr{F}, \{{\mathscr{F}}_t\},\P)}
\def\filqspace{(\Oo, \mathscr{H}, \{{\mathscr{H}}_t\},\Q)}
\def\filTpspace{(\Ot, \mathscr{F}_T, \{{\mathscr{F}}_t\},\Ptn)}

\def\RU{R_U}
\def\RV{R_V}
\def\RVi{R_{\Vi}}
\def\RVo{R_{\Vo}}
\def\RVt{R_{\Vt}}

\def\Uo{U^{(1)}}
\def\Vo{V^{(1)}}

\def\Ut{U^{(2)}}
\def\Vt{V^{(2)}}

\def\Ui{U^{(i)}}
\def\Vi{V^{(i)}}

\def\Wm{\mathscr W}

\def\WttB{\Wtm_t(B)}
\def\WtitB{\Wtim_t(B)}
\def\Wtmn{\Wtim_n}
\def\WtitBn{\Wtim_{t\wedge\taunU}(B)}
\def\solUWt{(U,\Wtmn)}

\def\Wtm{\Wm^{(2)}}
\def\Wtim{\tilde{\Wm}}
\def\Witim{\Wtim^{(i)}}
\def\solViWin{(\Vi,\Wmin)}

\def\solViWiti{(\Vi,{\Witim})}
\def\tauni{\tau_n^{(i)}}
\def\tauno{\tau_n^{(1)}}
\def\taunt{\tau_n^{(2)}}
\def\taunU{\tau_n^U}
\def\taunV{\tau_n^V}

\def\Wmi{\Wm^{(i)}}
\def\Wmit{\tilde{\Wm}^{(i)}}
\def\Wmin{\Wm_n^{(i)}}
\def\WitB{\Wm_t^{(i)}(B)}
\def\WitBn{\Wm_{t\wedge\tauni}^{(i)}(B)}
\def\WtitB{\tilde{\Wm}_t(B)}
\def\WittB{\tilde{\Wm}_t^{(i)}(B)}
\def\BL{{\mathscr B}(\S)}
\def\BR{{\mathscr B}(\S)}

\def\solViWiti{(\Vi,{\Witim})}
\def\RNbTtauniLYiWit{\Upsilon _{T\wedge\tauni}^{\RVi,\tilde {\Wm}^{(i)}}(\S)}
\def\RNbttauniBYiWit{\Upsilon _{t\wedge\tauni}^{\RVi,\tilde {\Wm}^{(i)}}(B)}

\def\RNfTtauniLYiWi{\Xi_{T\wedge\tauni}^{\RVi,{\Wm}^{(i)}}(\S)}

\def\RNfTtaunoLYoWo{\Xi _{T\wedge\tauno}^{\RVo,{\Wm}^{(1)}}(\S)}
\def\RNfTtauntLYtWt{\Xi _{T\wedge\taunt}^{\RVt,{\Wm}^{(2)}}(\S)}

\def\RNfTtaunW{\Xi_{T\wedge\taunU}^{\RU,{\Wm}}(\S)}
\def\RNfTtaunWt{\Xi_{T\wedge\taunU}^{\RU,{\Wtm}}(\S)}
\def\Xx{X^x}
\def\Xix{X^{-ix}}

\def\KSfBTP{{\mathbb {A}^{X,B}_{\un}(t,x)}}

\begin{document}
\title[From the L-KS kernel to noisy L-KS PDE{\scriptsize s} in multidimensions] {L-Kuramoto-Sivashinsky SPDE{\scriptsize s} in one-to-three dimensions: L-KS kernel, sharp H\"older regularity, and Swift-Hohenberg law equivalence}

\author{Hassan Allouba}
\address{Department of Mathematical Sciences, Kent State University, Kent,
Ohio 44242}
\email{allouba@math.kent.edu}
\begin{abstract} 
Generalizing the L-Kuramoto-Sivashinsky (L-KS) kernel from our earlier work, we give a novel explicit-kernel formulation useful for a large class of fourth order deterministic, stochastic, linear, and nonlinear PDEs in multispatial dimensions.  These include pattern formation equations like the Swift-Hohenberg and many other prominent and new PDEs.  We first establish existence, uniqueness, and sharp dimension-dependent spatio-temporal H\"older regularity for the canonical (zero drift) L-KS SPDE, driven by white noise on $\{\Rp\times\Rd\}_{d=1}^{3}$.  The spatio-temporal H\"older exponents are exactly the same as the striking ones we proved for our recently introduced Brownian-time Brownian motion (BTBM) stochastic integral equation, associated with time-fractional PDEs.  The challenge here is that, unlike the positive BTBM density, the L-KS kernel is the Gaussian average of a modified, highly oscillatory, and complex  Schr\"odinger propagator.  We use a combination of harmonic and delicate analysis to get the necessary estimates.   Second, attaching order parameters $\vepo$ to the L-KS spatial operator and $\vept$ to the noise term, we show that the dimension-dependent critical ratio $\vept/\vepo^{d/8}$ controls the limiting behavior of the L-KS SPDE, as $\vepo,\vept\searrow0$; and we compare this behavior to that of the less regular second order heat SPDEs.  Finally, we give a change-of-measure equivalence between the canonical L-KS SPDE and nonlinear L-KS SPDEs.  In particular, we prove uniqueness in law for the Swift-Hohenberg and the law equivalence---and hence the same H\"older regularity---of the Swift-Hohenberg SPDE and the canonical L-KS SPDE on compacts in one-to-three dimensions. 
\end{abstract}
\maketitle
\newpage
\tableofcontents
\section{Introduction and statements of results}
We give a novel, unifying, and very useful explicit-kernel (mild) formulation for a large class of linear, nonlinear, deterministic, and stochastic fourth order PDEs that includes many new, as well as prominent, equations.  We focus in this article on the L-Kuramoto-Sivashinsky (L-KS) stochastic PDEs\footnote{The name comes from the fundamental role of the linearized KS operator $-\tfrac\vep8\lpa\lap+2\vth\rpa^{2}$ in the nonlinear SPDE \eqref{nlks}.} (SPDEs):
\begin{equation} \label{nlks}
 \begin{cases} \displaystyle\frac{\partial U}{\partial t}=
-\tfrac\vep8\lpa\lap+2\vth\rpa^{2}U+b(U)+a(U)\frac{\partial^{d+1} W}{\partial t\partial x}, & (t,x)\in(0,+\infty )\times\Rd;
\cr U(0,x)=\unx, & x\in\Rd.
\end{cases}
\end{equation}
where $(\vep,\vth)\in(0,\infty)\times\R$ is a pair of parameters, $a,b:\cC\to\cC$ and $\un:\Rd\to\R$ are Borel measurable, and ${\partial^{d+1} W}/{\partial t\partial x}$ is the space-time white noise corresponding to the real-valued Brownian sheet\footnote{As in Walsh \cite{W}, we treat space-time white noise as a continuous orthogonal martingale measure, and we denote it by $\sW$.}  $W$ on $\Rp\times\Rd$, $d=1,2,3$.  In particular $b(u)$ may be a polynomial of (1) Allen-Cahn type $b(u)=\sum_{k=0}^{2p-1}c_{k}u^{k}$ for $p\in\N$ and for $c_{2p-1}<0$, to get many interesting fourth order SPDEs with an Allen-Cahn type nonlinearity (including a generalized Swift-Hohenberg SPDE when $\vep,\vth>0$), or of (2) KPP type $b(u)=\sum_{k=0}^{2p}c_{k}u^{k}$ for $p\in\N$ and for $c_{2p}<0$, to get new  fourth order SPDEs with a KPP type nonlinearity\footnote{The corresponding deterministic PDEs are, of course, obtained by simply setting $a\equiv0$.}.   We then use our explicit-kernel formulation to obtain, among other things, existence, uniqueness, and dimension-dependent H\"older regularity results with sharp spatio-temporal H\"older exponents for versions of the fourth order SPDE \eqref{nlks}.  More specifically, our first result \thmref{Holreg} establishes existence, uniqueness, and sharp dimension-dependent H\"older regularity for the zero drift ($b\equiv0$ or canonical L-KS SPDE) fixed $(\vep,\vth)$ version of \eqref{nlks}; \thmref{intense} gives dimension-dependent order parameters limiting results on the competing interaction between the linearized Kuramoto-Sivashinsky (L-KS) operator $-\tfrac\vep8\lpa\lap+2\vth\rpa^{2}$ and the white noise $\sW$ in \eqref{nlks} with $b\equiv0$ and $\vth$ fixed, and it gives the precise order parameters rate controlling whether the $L^{p}$ distance between an L-KS SPDE and its corresponding L-KS PDE uniformly vanishes (as the rate goes to zero) or whether there is a finite-time $L^{2}$ blowup of L-KS SPDEs (as the rate goes to infinity); and \thmref{com} adapts our earlier space-time change of measure results, with widely applicable conditions---from the second order SPDEs case \cite{Acom,Acom1,Acom2} to our fourth order SPDEs setting here---to transfer uniqueness in law and establish the law equivalence between the zero drift ($b\equiv0$) and the nonlinear nonzero drift versions of \eqref{nlks} on $\lbr\Rp\times\Rd\rbr_{d=1}^{3}$ and compact rectangles thereof.  This allows us to transfer almost sure properties of solutions--- including regularity---between linear and nonlinear L-KS SPDEs in spatial dimensions $d=1,2,3$.  An important special case covered by \thmref{com} is the aforementioned Swift-Hohenberg SPDE.  

We note here that the deterministic Swift-Hohenberg PDE (both real and complex) models numerous pattern formation phenomena in physics, chemistry, and optics (see e.g. \cite{BucLinPar,CroHoh,HohSwi,PelRot,StaSan,SwiHoh,Tri,DuWaSun}). These include the Taylor-Couette flow, the Rayleigh-B\'enard convection problem in a horizontal fluid layer in the gravitational field, large-scale flows and spiral core instabilities, and some chemical reactions.  Also, in optics, this equation is connected to spatial structures in large aspect lasers and synchronously pumped optical parametric oscillators.  The noisy Swift-Hohenberg PDE (or SPDE) treated here in \thmref{com} and \coref{SHcom} is at least as interesting and applicable.
We also remark that we use our kernel representational approach in separate papers to investigate time asymptotics and other qualitative behavior of a class of fourth order equations with different nonlinearities. 
\bnt\lbl{spdesnot}
We sometimes denote by $\eLKSabun$ the SPDE \eqref{nlks} on subsets of $\lbr\Rp\times\Rd\rbr_{d=1}^{3}$.  Similarly, the zero drift case is denoted by $\eLKSazun$. 
\ent
Before precisely stating our results, it is instructive to motivate and put together the building blocks---and give the different interesting links---in our approach; and then give our solution formulation and definition.  
\subsection{The L-KS PDE and L-KS kernel}\lbl{lkspdekersec}
We start with what we call the $(\vep,\vth)$ linearized Kuramoto-Sivashinsky (L-KS) PDE
\beq\lbl{vepvthlkspde} 
\frac{\partial u}{\partial t}=-\tfrac\vep8\lpa\lap+2\vth\rpa^{2}u; \ t>0, x\in\Rd, d\in\N=\{1,2,3\ldots\},
\eeq
  and we observe that it is a fundamental part of and/or intimately connected to a large family of interesting linear, nonlinear, deterministic, and stochastic PDEs.  This family includes, but is not limited to, both prominent and new compelling fourth order equations---including pattern formation equations---like   
\ben
\item the PDE\footnote{Throughout the article we alternate freely between the notations $\pa_{x}$ and $\pa/\pa x$ (or $d/dx$) for partial (or full) derivatives, with respect to any variable $x$, for aesthetic and typesetting reasons.} $\pa_{t}u=-\tfrac\vep8\lpa\lap+2\vth\rpa^{2}u+b(u)$, which includes the Swift-Hohenberg (SH) PDE (when $b(u)$ is an Allen-Cahn type nonlinearity and $\vep,\vth>0$) as well as many other interesting equations;
\item variants/versions of the Cahn-Hilliard PDE $\pa_{t}u=-\frac\vep8\lap^{2}u+\lap b(u)$, where $b$ may be an Allen-Cahn type nonlinearity $b(u)=\sum_{k=0}^{2p-1}c_{k}u^{k}$, $p\in\N$, and $c_{2p-1}<0$, $\vep>0$;
\item variants/versions of the Kuramoto-Sivashinsky (KS) PDE like 
\beq\lbl{ksandmodks}
\bsp
\partial_{t} u&=-\tfrac\vep8\lap^{2}u-\alpha_{1}\lap u-\alpha_{2}\nabla(u,u), \mbox{ and }
\\\partial_{t} u&=-\tfrac\vep8\lpa\lap+2\vth\rpa^{2}u-\frac12\nabla(u,u),
\end{split}
\eeq
where $\alpha_{1},\alpha_{2},\vep,\vth>0$;
\een 
and the stochastic versions of all the above PDEs, as well as many more new and intriguing fourth order (S)PDEs.  Some of these nonlinear equations mentioned have been studied, and continue to be studied, extensively in the deterministic literature (e.g., \cite{CT,CF,FK,FT,JRT,T,TW} and the SH references above) and is catching up on the still growing stochastic side (e.g., \cite{DaDeb,Duanebook,WaWaWa,DuWuCui,DuWaSun}), where the effect of the noise on the qualitative behavior of the underlying PDEs is of great interest.  When they are studied in the presence of a driving space-time white noise---with only few exceptions like \cite{DaDeb} and, recently, our work on higher order stochastic equations \cite{Atfhosie,Abtbmsie,Abtpspde}---these fourth order equations are invariably restricted to one spatial dimension $d=1$.   On the other hand, in our earlier work \cite{Abtp1,Abtp2,Aks}; we introduced and connected a large class of processes---in which the time parameter is replaced in different ways by a Brownian motion---to new memory-preserving (\emph{memoryful}) fourth order PDEs and to the linearized KS PDE \eqref{vepvthlkspde} with $\vep=\vth=1$:
\begin{equation} \label{lkspde}
 \begin{cases} {\ds\frac{\partial u}{\partial t}}=
-\frac18\lap^2u-\frac12\lap u-\frac12u, & (t,x)\in(0,+\infty )\times\Rd;
\cr u(0,x)=\unx, & x\in\Rd,
\end{cases}
\end{equation}
{\it in all spatial dimensions\footnote{This is important to note since one of the major challenges in the study of the nonlinear KS equation is that the existence of solutions in spatial dimensions $d\ge2$ is unsettled, even in the noiseless deterministic case (see \cite{T}).} $d\ge1$} for suitably regular initial data $\un$.  At the heart of our approach in \cite{Aks} is the kernel $\KKStxy$---associated with what in \cite{Aks} we call the imaginary-Brownian-time-Brownian-angle process (IBTBAP)---defined by
\beq\lbl{ibtbapkernel}
\bc
\Mpisxy\eqdef{\exp\lpa\i s\rpa}\df{\e^{-|x-y|^2/2\i s}}{{\lpa2\pi \i s \rpa}^{d/2}},\cr
\KKStxy\eqdef\ds\int_{-\infty}^0\Mpisxy\ptsz ds+\int_{0}^\infty\Mpisxy\ptsz ds
\ec
\eeq
where $\i=\sqrt{-1}$ and $\ptsz=\tf{\e^{-s^2/2t}}{\sqrt{2\pi t}}$.  Since $\KKStx$, obtained by setting $y=0\in\Rd$ in $\KKStxy$, is the fundamental solution of the L-KS PDE in \eqref{lkspde}, we also call it the L-KS kernel\footnote{See also \secref{kerforsec} and Section 2 below for a simpler form of $\KKStx$ and its connection to \eqref{lkspde} via Fourier transforms.}.  
Quantum mechanics experts will quickly recognize that, except for the $\exp{(\i s)}$ angle term, $\Mpisxy$ in the definition of the L-KS kernel in \eqref{ibtbapkernel} is a $d$-dimensional version of the free propagator associated with Schr\"odinger's equation.  It is then proved in Theorem 1.1 of \cite{Aks} that, for\footnote{The compact support assumption on $\un$ here and in \thmref{Holreg} below is for convenience only and may be replaced with more relaxed integrability conditions \`a la those given for the Brownian-time Brownian sheet in \cite{AN}.} $\un\in\C^{2,\gamma}_c(\Rd;\R)$, 
\beq\lbl{detlksint}
u(t,x)=\int_{\Rd}\KKStxy\uny dy
\eeq
solves the linearized Kuramoto-Sivashinsky PDE \eqref{lkspde} and hence that the kernel $\KKStx$ solves the PDE L-KS in \eqref{lkspde} with initial condition $\delta(x)$. 
\subsection{Imaginary-Brownian-time-Brownian-angle and Schr\"odinger links}
An important intuitive ingredient in the formulation and proof of Theorem 1.1 of \cite{Aks}, and in arriving at the kernel $\KKStx$, was the use of the intimate connection between the Brownian-time processes and their densities in \cite{Abtp1,Abtp2} and the imaginary-Brownian-time-Brownian-angle process and its kernel\footnote{In particular, Theorem 1.2 in \cite{Abtp2} was crucial in arriving at our IBTBAP and its kernel $\KKStx$.  Of course, the L-KS kernel $\KKStx$ is \emph{not} a proper probability density in the standard sense.  But, it has a nice Fourier transform, as we shall see shortly.} in \cite{Aks}.
Our IBTBAP process, starting at $\un:\Rd\to\R$, was given in \cite{Aks} by
\begin{equation}
\KSfBTP\eqdef\begin{cases} \un\lpa X^{x}(\i B(t))\rpa \exp\left(\i B(t)\right),  & B(t)\ge0;\cr
\un\lpa\i X^{-\i x}(-\i B(t))\rpa\exp\left(\i B(t)\right),  & B(t)<0;
\end{cases}
\label{KSprocess}
\end{equation}
where the process $X^{x}$ is an $\Rd$-valued Brownian motion (BM) starting from $x\in\Rd$, $\Xix$ is an independent $i\Rd$-valued BM starting at $-ix$ (so that $i\Xix$ starts
at $x$), and both are independent of the inner standard $\R$-valued Brownian motion $B$ starting from $0$.    The clock of the outer Brownian motions $\Xx$ and $\Xix$ is replaced
by an imaginary positive Brownian time; and the angle of $\KSfBTP$ is the Brownian motion $B$.
We think of the imaginary-time processes $\{\Xx(is),s\ge0\}$ and $\{i\Xix(-is),s\le0\}$ as having the same complex Gaussian distribution on $\Rd$ with the corresponding complex distributional density (or Schr\"odinger propagator)
$${K}^{\text{\tiny{\sc{SP}}}^d}_{\i s;x,y}=\frac{1}{({2\pi i s})^{d/2}}e^{-|x-y|^2/2is}.$$
We may then think of $u$ in \eqref{detlksint} in terms of complex expectation by first conditioning on $B(t)=s$ and then removing the conditioning (by integrating over $s$) and defining
$u(t,x)\eqdef\E^{\cC}\left[\KSfBTP\right]$.  Viewed this way, $\KKStx$ is the expectation kernel of the IBTBAP.    Since $\Mpisxy=\e^{\i s}{K}^{\text{\tiny{\sc{SP}}}^d}_{\i s;x,y}$ is obtained by giving the propagator ${K}^{\text{\tiny{\sc{SP}}}^d}_{\i s;x,y}$ an extra angle $s\in\R\setminus\{0\}$, where $s$ is also the real-valued time on the imaginary axis ($\i s$), we call $\Mpisxy$ the $d$-dimensional $\R$-time-angled propagator.  The L-KS kernel $\KKStx$ in \eqref{ibtbapkernel} is thus the Gaussian average of an $\R$-time-angled Schr\"odinger propagator\footnote{In our fourth order setting we have two notions of time: the standard time $t$ and the Brownian-time and Brownian-angle $B(t)$ (each $s$ in $\Mpisxy$ represents a possible value for the BM $B$ in our IBTBAP at some time $t$, $B(t)=s\in\R$).}.  
\subsection{The $(\vep,\vth)$ L-KS kernel formulation}\lbl{kerforsec}
In this article, we start by using our L-KS kernel to formulate the notion of a mild kernel solution to the $(\vep,\vth)$ L-KS (S)PDEs in \eqref{nlks}.  We first generalize slightly $\KKStx$ in \eqref{ibtbapkernel} by scaling the time $t$ with a parameter $\vep>0$ and scaling the angle $s$ in the $\R$-time-angled propagator $\Mpisxy$ by another parameter $\vth\in\R$ to obtain the $(\vep,\vth)$ L-KS kernel\footnote{Clearly, using the notation of the-just-introduced $(\vep,\vth)$ L-KS kernel, we note that $\KKStx=\KKSepthootx$.}

\beq\lbl{vepvthLKS}       
\KKSepthtxy=\int_{-\infty}^{0}\df{\e^{i\vth s} \e^{-|x-y|^2/2\i s}}{{\lpa2\pi \i s \rpa}^{d/2}}\peptsz ds+\int_{0}^{\infty}\df{\e^{i\vth s} \e^{-|x-y|^2/2\i s}}{{\lpa2\pi \i s \rpa}^{d/2}}\peptsz ds,
\eeq 
which, when setting $y=0\in\Rd$, is the fundamental solution $\KKSepthtx$ to the $(\vep,\vth)$ L-KS PDE in equation \eqref{vepvthlkspde}\footnote{See Section 2 for a Fourier argument.  We also briefly note that with $\vep>0$, we always have the dissipative negative biLaplacian $-\Delta^{2}$ in \eqref{vepvthlkspde}.  On the other hand, the case $\vth<0$ leads to a dissipative second order $\Delta$; whereas $\vth>0$ leads to the non-dissipative second order $-\Delta$, which is the case in L-KS PDEs like the Swift-Hohenberg and Kuramoto-Sivashinsky for example.}.   As we will see in \secref{Foursec}, despite the involved expression in \eqref{vepvthLKS}, the kernel $\KKSepthtx$ has a rather nice (and revealing) Fourier transform:
\beq\lbl{FTkerintro}
\FKKSepthtxi=\lpa2\pi\rpa^{-\frac d2}\e^{-\frac{\vep t}{8}\left( -2\vth+\lab\xi\rab^{2} \right) ^{2}};\ \vep>0,\ \vth\in\R,
\eeq
which, upon inverting yields the simpler form of $\KKSepthtx$
\beq\lbl{iFTintro}
\bsp
\KKSepthtx&=(2\pi)^{-d}\int_{\Rd}\e^{-\frac{\vep t}{8}\left( -2\vth+\lab\xi\rab^{2} \right) ^{2}}\e^{\i\xi\cdot x} d\xi; 
\\&= (2\pi)^{-d}\int_{\Rd}\e^{-\frac{\vep t}{8}\left( -2\vth+\lab\xi\rab^{2} \right) ^{2}}\cos\lpa{\xi\cdot x}\rpa d\xi;\ \vep>0,\ \vth\in\R.
\end{split}
\eeq  
The last equality in \eqref{iFTintro} follows since $\int_{\Rd}\e^{-\frac{\vep t}{8}\left( -2\vth+\lab\xi\rab^{2} \right) ^{2}}\sin\lpa{\xi\cdot x}\rpa d\xi=0$.  Thus, the effect of the Gaussian average of the propagator in \eqref{vepvthLKS} is to ``average out'' the imaginary part of the kernel, and $\KKSepthtx$ is real-valued.   We now give the new kernel formulation for the class of $(\vep,\vth)$ L-KS (S)PDEs \eqref{nlks} which includes, among many other (S)PDEs, the Swift-Hohenberg (S)PDE.
\bdf[$(\vep,\vth)$ L-KS kernel (mild) formulation of $(\vep,\vth)$ L-KS (S)PDEs \eqref{nlks}]\lbl{kskersoldef}
Fix $\vep>0$ and $\vth\in\R$.  We call the pair $(U,\sW)$ on a usual probability space\footnote{We assume throughout the article that filtrations $\{\sFt\}_{t\ge0}$ satisfies the usual conditions, and we often simply say that $U$ is a kernel solution to \eqref{nlks} to mean the same as the definition of a mild solution above.} $\OFFtP$ a $(\vep,\vth)$ L-KS kernel solution to \eqref{nlks} on $\Rp\times\Rd$ whenever $\sW$ is a space-time 
white noise on $\Rp\times\Rd$; the random field $U$ is progressively measurable, and with $U(0,x)=\unx$; and the pair $(U,\sW)$ satisfies the $(\vep,\vth)$ L-KS kernel (mild) formulation:
\beq\lbl{ibtbapsol}
\bsp
U(t,x)=&\intrd\KKSepthtxy \uny dy
\\&+ \intrdzt \KKSepthtsxy\lbk b(\usy)ds dy +a(\usy)\sW(ds\times dy)\rbk
\end{split}
\eeq
for $t>0$ and for every (or almost every) $x\in\Rd$, almost surely $\P$.  Weak and strong---in the probability sense---solutions are defined in the usual way:  we call a $(\vep,\vth)$ L-KS kernel solution weak if the white noise $\sW$ and $\OFFtP$ on which it's defined are freely chosen---along with $U$---so as to satisfy \eqref{ibtbapsol}; and the solution is strong if $\sW$ and $\OFFtP$ are fixed and $\{\sFt\}$ is the augmentation of the natural filtration for $\sW$ under $\P$.  The solution is continuous if $U$ has continuous paths on $\Rp\times\Rd$ almost surely $\P$. 

Uniqueness in law holds for \eqref{nlks} if the laws\footnote{All solutions $U$ in this article have continuos paths ($U\in\C(\dD;\R)$, where $\dD\subset\Rp\times\Rd$).  The law $\LawUP$ of the random field $U$ under $\P$ is the probability measure induced on the Borel $\sigma$-field of continuous function by the recipe: $\LawUP(\Lambda)=\P[U\in\Lambda]$, $\Lambda\in\sB(\C(\dD;\R))$  }  $\LawUiPi$ of $\Ui$
under $\Pii$; $i=1,2$, are the same whenever
$(\Ui,\sWi)$, $\OFFtPi$; $i=1,2$, are $(\vep,\vth)$ L-KS kernel solutions to \eqref{nlks}.  Pathwise uniqueness holds if $\Uo$ and $\Ut$ are $\P$-indistinguishable ($\P\lbk\Uo=\Ut\rbk=1$) whenever $(\Ui,\sW)$ are $(\vep,\vth)$ L-KS kernel solutions to \eqref{nlks} with respect to the same white noise $\sW$ and on the same probability space $\OFFtP$.
 
A $(\vep,\vth)$ L-KS kernel solution $U$ to the deterministic version of \eqref{nlks} is obtained from \eqref{ibtbapsol} by setting $a\equiv0$.  \edf
\brm\lbl{kernlwideapp}
Although we focus in this article on the SPDE in \eqref{nlks} on $\{\Rp\times\Rd\}_{d=1}^{3}$ and subsets thereof, the utility of our new L-KS kernel formulation goes well beyond just \eqref{nlks}.   We show separately how to adapt it to formulate the class of PDEs discussed in \secref{lkspdekersec} (Cahn-Hilliard, Kuramoto-Sivashinsky, and many other fourth order equations) and their stochastic versions.   We also use this explicit kernel approach in separate articles to analyze the time-asymptotic\footnote{See \cite{AL04} for a similar approach in studying random attractors for the second order Allen-Cahn case.} and other qualitative behaviors of several fourth order L-KS type equations.
\erm
\subsection{Three main theorems and a Swift-Hohenberg corollary}
In this article, we establish three main theorems on versions of the $(\vep,\vth)$ L-KS SPDE \eqref{nlks}.  We now detail and state our main results.
\subsubsection{\thmref{Holreg}: existence, uniqueness, and sharp dimension-dependent H\"older regularity}  In our first result;  we obtain sharp, dimension-dependent, spatio-temporal H\"older continuity regularity results for the L-KS SPDE \eqref{nlks} with $b\equiv0$ (zero drift or canonical L-KS SPDE):
\begin{equation} \label{nlkszd}
 \begin{cases} \displaystyle\frac{\partial U}{\partial t}=
-\tfrac\vep8\lpa\lap+2\vth\rpa^{2}U+a(U)\frac{\partial^{d+1} W}{\partial t\partial x}, & (t,x)\in(0,+\infty )\times\Rd;
\cr U(0,x)=\unx, & x\in\Rd.
\end{cases}
\end{equation}

In particular, for any fixed $\vep,T>0$ and $\vth\in\R$, we obtain the existence of a unique real-valued solution $U$ that is $L^p(\Omega)$-bounded on $[0,T]\times\Rd$ for all $p\ge2$ and that has H\"older continuous paths in time and space.  In time, the H\"older exponent is $\gamma_{t}\in\lpa0,(4-d)/8\rpa$ and in space it is $\gamma_{s}\in\lpa0,[(4-d)/2]\wedge1\rpa$, for spatial dimensions $d=1,2,3$.   We first obtained the same striking spatio-temporal H\"older regularity profile in \cite{Abtbmsie} for a different class of memoryful fourth order stochastic integral equations (SIEs) associated with the Brownian-time Brownian motion (BTBM)---see \cite{Abtp1,Abtp2,Abtpspde} and the discussion in \cite{Abtbmsie}---which we introduced as BTBM SIEs.  What is remarkable about this H\"older regularity profile is that, not only random field solutions exist in spatial dimensions $d=1,2,3$ (not just for $d=1$) in the presence of the rough driving space-time white noise\footnote{This is in contrast to second order PDEs driven by space-time white noise whose random field solution exists only in $d=1$.  Also, it is noteworthy that, with very few exceptions (e.g., \cite{Atfhosie,Abtbmsie,Abtpspde,DaDeb}), space-time white noise driven SPDEs, even higher order ones,  are not treated in more than one spatial dimension.}, but these random field solutions are spatially twice as smooth as the underlying Brownian sheet\footnote{Our article \cite{Abtbmsie} gave the first example of space-time white noise driven equations whose solutions are smoother in either time or space than the underlying Brownian sheet corresponding to the driving white noise.} in $d=1,2$.  In the followup article \cite{Atfhosie}, we showed that this third dimensionality limit on random field existence and the above spatial H\"older smoothness are maximal\footnote{Maximal if the spatial dimension is integer.} in equations driven by a space-time white noise that are first order in time and high order in space---no matter how high the order is in these equations.   

Although our L-KS SPDEs here have the same spatio-temporal H\"older profile as the BTBM SIE of \cite{Abtbmsie}, proving it by directly adapting our methods in \cite{Abtbmsie} to the L-KS kernel is demanding.  The difficulty lies in the fact that the L-KS kernel in \eqref{vepvthLKS} is the Gaussian average of the  highly oscillatory angled complex propagator; whereas the BTBM probability density 
\beq\lbl{btpden}\KBtxy=2\int_0^\infty \df{e^{-|x-y|^2/2s}}{(2\pi s)^{{d}/{2}}}\df{e^{-s^2/2t}}{\sqrt{2\pi t}} ds
\eeq
is a Gaussian average of another non-oscillatory Gaussian density.   Also, the L-KS kernel is not a proper probability density function as the BTBM density.  Thus, we proceed differently by applying a harmonic analytic step to the L-KS kernel at the outset.  This turns out to be a useful first step towards obtaining the required regularity estimates.  We then use delicate analysis, including comparing the nonzero $\vth$ angle case to that of the simpler $\vth=0$ case and adapting the probabilistic-analytic arguments from \cite{Abtbmsie} to our setting here, to prove the estimates needed for the proof of \thmref{Holreg}.
 
In the process, we give a harmonic analytic explanation of why the BTBM density---which is associated with the quite different memory-preserving \emph{positive biLaplacian} fourth order PDE
\beq\lbl{btppdedet}
\bc
\df{\pa u}{\pa t}= \df{\Delta\un }{\sqrt{8\pi t}}+\df18\Delta^2 u;&(t,x)\in(0,\infty)\times\Rd, \cr
u(0,x) = \unx;&x\in\Rd,
\ec
\eeq
and its equivalent\footnote{The connection of BTBM to fourth order PDEs (including \eqref{btppdedet}) was first established in our papers \cite{Abtp1,Abtp2}.  Also, their connection to time-fractional PDEs was first established implicitly via the BTBM connection to the half derivative generator in \cite[Theorem 0.5]{Abtp1}.  The equivalence between a large class of time-fractional and higher order PDEs with memory, including the equivalence between \eqref{btppdedet} and \eqref{BTBM-half} when $\un\in\C_{b}^{2,\gamma}$, was established explicitly in \cite[Theorem 3.1--Theorem 3.6]{BMN} (see also \cite{BOap09} and \cite{AN}).    For a further discussion of interesting aspects of these PDEs and their history see also \cite{Abtp1} and the introduction in \cite{Atfhosie}.  In the recent multiparameter-time case the reader is referred to  \cite{Abtbs,AN}.  The BTBM scaling and its nonstandard PDEs  connection have now attracted a lot of attention, even outside probability and PDEs, as evidenced by the recent physics and mathematical finance articles \cite{CalEicSau,Carr}.} time-fractional PDE
\begin{equation}\label{BTBM-half}
\begin{cases}
\displaystyle\pa_{t}^{\frac12} u=\frac{1}{\sqrt{8}}\lap u;& t\in(0,\infty),x\in\Rd,\\
u(0,x)=\unx;&x\in\Rd,
 \end{cases}
\end{equation} 
where $\pa_{t}^{\frac12}$ is the Caputo fractional derivative---has the same regularizing effect as that of the L-KS kernel and its associated PDE \eqref{lkspde}.   
This harmonic explanation is given in Section 2 below.  A different probabilistic heuristic argument was given in \cite{Abtbmsie} as to why the BTBM-SIEs in \cite{Abtpspde,Abtbmsie} are cousins of the L-Kuramoto-Sivashinsky SPDEs here.  The Fourier transform of the L-KS kernel, and its inverse, are also used in Section 2 to sketch a different proof of the L-KS PDE \eqref{lkspde} connection to the L-KS kernel, first proved differently in \cite{Aks}.
 
The regularity---and other qualitative behavior---results carry over to a large class of nonlinear L-KS SPDEs like \eqref{nlks} and others intimately connected to the linear KS operator $-\tfrac\vep8\lpa\lap+2\vth\rpa^{2}$.  Some of these are illustrated in \thmref{com}, which is possible by adapting our earlier change of measure results \cite{Acom,Acom1,Acom2} from the second order to the fourth order settings.  Understanding the L-KS PDE \eqref{lkspde} and SPDE \eqref{nlkszd} is thus very useful in understanding a large class of nonlinear L-Kuramoto-Sivashinsky equations; including the Swift-Hohenberg and its generalization \eqref{nlks}, variants of the KS, and many more.  

Throughout this paper, fix an arbitrary $T>0$, and let $\T=[0,T]$.  We denote by $\H^{\gamma_t^{-},\gamma_s^{-}}(\T\times\Rd;\cC)$ the space of real-valued locally H\"older functions on $\T\times\Rd$ whose time and space H\"older exponents are in $(0,\gamma_t)$ and $(0,\gamma_s)$, respectively.  We now state our first existence, uniqueness, and regularity result\footnote{We remind the reader again that the compact support condition on $\un $ may be      replaced with more relaxed integrability conditions like those given for the Brownian-time Brownian sheet in \cite{AN}.}.
\bfr
\bthm[Existence/uniqueness and sharp H\"older regularity for the canonical $(\vep,\vth)$ L-KS \eqref{nlkszd}  in $d=1,2,3$]\lbl{Holreg}
 Fix $\vep>0$ and $\vth\in\R$.  Assume that
\beq\lbl{lcnd}\tag{\mbox{Lip}}
\bc
(a) &\lab a(u)-a(v)\rab\le C\lab u-v\rab \mbox \quad u,v\in\cC;\\
(b) &a^2(u)\le C (1+\lab u\rab^2),\quad u\in\cC;\\
(c) &u_0\in \mathrm{C}_{c}^{2,\gamma}(\Rd;\R)\mbox{ and nonrandom } \forall\ d\in\{1,2,3\}.
\ec
\eeq
Then, there exists a strong $(\vep,\vth)$ L-KS kernel solution solution $(U,\sW)$ to the L-KS SPDE \eqref{nlkszd} on $\Rp\times\Rd$, for $d=1,2,3$, which is $L^p(\Omega)$-bounded on $\T\times\Rd$ for all $p\ge2$ $(M_{p}(t):=\sup_{x\in\Rd}\E\lab U(t,x)\rab^{p}\le C_{T}$ for $t\in\T)$ and 
$$U\in\H^{{\tf{4-d}{8}}^{-},{\lpa\tf{4-d}{2}\wedge 1\rpa}^{-}}\lpa\T\times\Rd;\cC\rpa;\mbox{ for  $d=1,2,3$, almost surely.}$$    
If $(V,\sW)$ is another such solution, with respect to the same white noise $\sW$, then, for any $d\in\{1,2,3\}$, $U$ and $V$ are indistinguishable: $\P\lbk U(t,x)=V(t,x); (t,x)\in\Rp\times\Rd\rbk=1.$  
\ethm
\efr
\brm\lbl{noliprem}
\ben\rencomrom
\item We note here that we can adapt our lattice arguments and K-martingale approach in \cite{Abtbmsie} to prove existence of lattice-limits solutions to our L-KS SPDE with the same $L^{p}$ and H\"older regularity as those of \thmref{Holreg} under the weaker non-Lipschitz conditions\footnote{These lattice arguments have their roots in our second order SPDE works \cite{Asdde1,Asdde2}.}  
 \beq\lbl{nlcnd}\tag{\mbox{NLip}}
\bc
\hspace{-2.5mm}&(a)\ a(u)\mbox{ is continuous in }u \mbox \quad u\in\cC;\\
\hspace{-2.5mm}&(b) \mbox{ and } (c)\mbox{ same as in \eqref{lcnd}};
\ec
\eeq
The details are left to the interested reader.  
\item If the Lipschitz condition on $a$ is only local, then the uniqueness assertion of \thmref{Holreg} holds for the L-KS SPDE \eqref{nlkszd}.  However, a local Lipschitz condition is not sufficient to guarantee global existence; and the existence of the unique solution $U$ and its H\"older regularity of \thmref{Holreg} hold only up to a random (possibly finite) dimension-dependent blowup time $\tau_{d}$, for which\footnote{Our work in a separate article indicates that, when $a(u)=|u|^{\alpha}$, there is a dimension-dependent critical blowup exponent $\alpha_{c}^{(d)}>1$ 
for $d=1,2,3$, such that
\ben\rencomrom
 \item if $\alpha>\alpha_{c}^{(d)}$, then the LKS SPDE \eqref{nlkszd} with diffusion coefficient $a(u)=|u|^{\alpha}$ blows up in finite time ($\tau_{d}<\infty$ in \eqref{blowup}); and 
 \item if $\alpha<\alpha_{c}^{(d)}$ then finite time blowup does not occur ($\tau_{d}=\infty$ in \eqref{blowup}).  \een
 More on that in the aforementioned upcoming article.}
\beq\lbl{blowup}
\P\lbk\lim_{t\nearrow\tau_{d}}\sup_{x}\lab U(t,x)\rab=\infty\rbk>0;\quad  d=1,2,3.
\eeq  
\item The H\"older exponents confirms our assertion in \cite{Abtpspde,Abtbmsie} about the intimate relation between our BTBM SIE there and the L-KS SPDE here.  They also tell us that the L-KS SPDE is much smoother than the second order heat SPDE whose solution $U\in\H^{{\tf{1}{4}}^{-},{\tf{1}{2}}^{-}}\lpa\T\times\R;\cC\rpa$ for only $d=1$.
\een
\erm
\subsubsection{\thmref{intense}: L-KS vs white noise in the limit, the case of vanishing intensities}  Consider the $(\vep,\vth)$ L-KS SPDE \eqref{nlkszd}. Fix $\vth\in\R$, let $\vep=\vepo>0$ be an order parameter, and attach another order parameter $\vept>0$ to the white noise term to obtain the L-KS SPDE
\begin{equation} \label{lksop}
 \begin{cases} \displaystyle\frac{\partial U}{\partial t}=
-\tfrac\vepo8\lpa\lap+2\vth\rpa^{2}U+\vept a(U)\frac{\partial^{d+1} W}{\partial t\partial x}, &\hspace{-2mm} (t,x)\in(0,+\infty )\times\Rd;
\cr U(0,x)=\unx, &\hspace{-2mm} x\in\Rd.
\end{cases}
\end{equation}
By changing the values of $\vepo$ and $\vept$, we are simply changing the intensities of the smoothing L-KS spatial operator $-\tfrac18\lpa\lap+2\vth\rpa^{2}$ and the roughening noise term, respectively.  \thmref{Holreg} gives us existense, uniqueness, sharp H\"older regularity, and $L^p(\Omega)$-boundedness on $[0,T]\times\Rd$ for all $T>0,p\ge2$, and hence no blowup, for \emph{any} fixed values of $\vepo,\vept$ on \emph{any} time interval $[0,T]$, under Lipschitz conditions on $a$.  Our second result, \thmref{intense}, is concerned with the interesting questions: what happens in the limit to the solution of \eqref{lksop} as $\vepo,\vept\searrow0$?  and how fast can $\vepo\searrow0$ relative to the speed at which $\vept\searrow0$ before the noise term dominates the L-KS operator?  \thmref{intense} answers these questions by giving a precise critical ratio, $\vept/\vepo^{d/8}$ for dimensions $d=1,2,3$, which controls the limiting behavior of \eqref{lksop} as $\vepo,\vept\searrow0$.   According to \thmref{intense}, this critical ratio tells us that if we let $\vepo,\vept\searrow0$ we should be careful \emph{not} to let $\vepo\searrow0$ so fast that $\vept/\vepo^{d/8}\nearrow\infty$ or the supremum---over space and time---of the solution's $L^{2}(\Omega)$ norm will grow to infinity on any time interval $[0,T]$.  On the other hand, this critical ratio also says that when we let $\vepo,\vept\searrow0$ so that $\vepo\searrow0$ slow enough that $\vept/\vepo^{d/8}\searrow0$, then the solution to \eqref{lksop} converges to its deterministic version ($a\equiv0$) $\vepo\searrow0$ limit in $L^{2q}(\Omega)$ for all $(q \ge 1)$.  The same phenomenon happens in the case of the less regular second order heat SPDE (see \remref{noliprem} (iii))
\begin{equation} \label{heat}
 \begin{cases} \displaystyle\frac{\partial U}{\partial t}=
-\tfrac\vepo2\lap U+\vept a(U)\frac{\partial^{2} W}{\partial t\partial x}, &\hspace{-2mm} (t,x)\in(0,+\infty )\times\R;
\cr U(0,x)=\unx, &\hspace{-2mm} x\in\R.
\end{cases}
\end{equation}
with a less forgiving critical ratio of $\vept/\vepo^{1/4}$ and only for $d=1$\footnote{As we remarked before, unlike our L-KS SPDEs here, random field solutions to these second order SPDEs are restricted to one dimensional space.}---as we showed in \cite[Theorem 1.8]{Asdde2}.   Comparing the two critical ratios for the L-KS SPDE \eqref{lksop} and the heat SPDE \eqref{heat}, we see another manifestation of the fact that the L-KS fourth order operator has  a substantially stronger regularizing effect on white noise---we can afford to turn it off much faster---than does the usual Laplacian in second order noisy heat PDEs.   For the critical ratio ${\vept}/{{ \vepo}^{1/8}}$ to grow to $\infty$ as $\vepo,\vept\searrow0$ in the one-dimensional L-KS case \eqref{lksop}, we have to ``turn off'' the smoothing L-KS operator (let $\vepo\searrow0$) at the fast rate of $\vepo=\vept^{8}z(\vept)$ for some function $z>0$ satisfying  $\lim_{\vept\searrow0}z(\vept)=0$ (e.g., $z(\vept)=\vept^{\delta}$ for $\delta>0$); whereas, in the second order heat SPDE case \eqref{heat}, the critical ratio there, ${\vept}/{{ \vepo}^{1/4}}\nearrow\infty$ as $\vepo,\vept\searrow0$ even if we turn off the smoothing Laplacian at the much slower rate of $\vepo=\vept^{4}z(\vept)$.  Also, for the critical ratio ${\vept}/{{ \vepo}^{1/8}}$ to converge to $0$ as $\vepo,\vept\searrow0$ in the one-dimensional L-KS case \eqref{lksop}, we can turn off the L-KS operator at a rate as high as $\vepo=\vept^{8}z^{-1}(\vept)$; however, in the second order heat SPDE case, we have to turn off the Laplacian at a much slower rate of no more than $\vepo=\vept^{4}z^{-1}(\vept)$ to cause the critical ratio there, ${\vept}/{{ \vepo}^{1/4}}$ to converge to $0$.  In $d=2,3$, while the heat SPDE \eqref{heat} does not have random field solutions, the L-KS SPDE \eqref{lksop} has H\"older continuous solutions (\thmref{Holreg}) and a critical ratio ${\vept}/{{ \vepo}^{d/8}}$.  In particular, in $d=2$, the critical ratio for the L-KS SPDE is ${\vept}/{{ \vepo}^{1/4}}$, the same as the one for the heat SPDE in $d=1$.  This is consistent with the fact that the effective H\"older exponent\footnote{The minimum of the themporal and spatial H\"older exponents.} of the L-KS SPDE in $d=2$ is the same as that for the heat SPDE in $d=1$, $(1/4)^{-}$.  In $d=3$, the critical ratio for the L-KS SPDE is ${\vept}/{{ \vepo}^{3/8}}$.  \thmref{intense} rigorizes the above heuristic for the L-KS SPDE  \eqref{lksop} and gives a precise meaning of ``too fast''  or ``too slow'' for the speed of the limit $\vepo\searrow0$, relative to the speed at which the intensity of the noise term $\vept$ is vanishing, by giving the precise critical ratio of $\vept$ to $\vepo$ that determines the limiting behavior as $\vepo,\vept\searrow0$.
\ig
In particular, \thmref{intense} below states that the $L^{2q}$ distance between the solution to the SPDE \eqref{lksop} and the solution to its deterministic version ($a\equiv0$) goes to zero, uniformly on $[0,T]\times\Rd$, as $\vepo,\vept, \mbox{ and }{\vept}/{{ \vepo}^{d/8}}\searrow0$ for $d=1,2,3$ and $q>1$.  It also says that the supremum---over time and space---of the $L^{2}(\Omega)$ norm of the solution grows to $\infty$ as $\vepo,\vept\searrow0$ such that the ratio ${\vept}/{{ \vepo}^{d/8}}\nearrow\infty$ for $d=1,2,3.$    

We now state our result which characterizes, through the critical ratio ${\vept}/{{ \vepo}^{d/8}}$, the limiting competition in the L-KS SPDE \eqref{lksop} between the regularizing effect of the spatial fourth order operator $-\tfrac18\lpa\lap+2\rpa^{2}$ as it pushes against the roughening effect of the space-time white noise. 
\fi
\ig
to study the  In \thmref{intense} we show that this competition is controlled in the limit by the critical ratio $\vept/\vepo^{d/8}$ as $\vepo,\vept\searrow0$ (or $\vepo,\vept\nearrow\infty$), for dimensions $d=1,2,3$.  
\fi 
  \bfr
\bthm[L-KS vs white noise in \eqref{lksop}: the critical order parameter ratio in $d=1,2,3$]\lbl{intense}
Fix $\vth\in\R$.  Assume that the conditions in \eqref{lcnd} are in force and that $(U_{\epsilon_1,\epsilon_2},\sW)$ is the unique strong solution to the L-KS SPDE \eqref{lksop}.
\ben\rencomrom
\item (Uniformly vanishing $L^{2q}$ distance between SPDE and PDE as ${\vept}/{{ \vepo}^{d/8}}\searrow0$) Suppose that $u_{\vepo}$ is the solution to the deterministic L-KS PDE obtained from \eqref{lksop} by setting $a\equiv0$, then 
$$\sup_{ 0\le t\le T}\sup_{x\in\Rd}\E\left|U_{\vepo,\vept}(t,x)-u_{\vepo}(t,x)\right|^{2q}\longrightarrow0; \forall  q\ge1, T>0$$
as $\vepo,\vept\searrow0, \mbox{ and }{\vept}/{{ \vepo}^{d/8}}\searrow0$ for $d=1,2,3$.
\item  (Supremum $L^{2}$ growth to infinity as ${\vept}/{{ \vepo}^{d/8}}\nearrow\infty$) Suppose there are constants $K_l,K_u>0$ such that $K_l\le a(v)\le K_u$
for all $v\in\mathbb{R};$   then,  
$$\sup_{0\le s\le T}\sup_{{x\in\Rd}}\mathbb{E} \lab U_{\epsilon_1,\epsilon_2}(s,x)\rab^{2}\nearrow\infty,\mbox{ $\forall T>0$.}$$ 
as $\vepo,\vept\searrow0$ such that the ratio ${\vept}/{{ \vepo}^{d/8}}\nearrow\infty$ for $d=1,2,3.$
\een
\ethm
\efr
Both \thmref{intense} above and \cite[Theorem 1.8]{Asdde2} capture the same common phenomenon\footnote{This phenomenon is also shared with the time-fractional/high order stochastic integral equations of \cite{Abtbmsie,Atfhosie} with the same (in \cite{Abtbmsie}) or more forgiving (in \cite{Atfhosie}) critical ratios as the L-KS, depending on the $\beta$ value.} for \eqref{lksop} and \eqref{heat}, respectively, associated with the vanishing of $\vepo$ and $\vept$. Each theorem gives a different critical ratio that is at the edge of the two extreme behaviors in \thmref{intense} (i) and (ii) and in \cite[Theorem 1.8]{Asdde2}.
\ig
which are heuristically understood as follows:    
\ben\rencomrom 
\item (subcritical case) if $\vepo\searrow0$ ``too slow'' compared to the speed at which $\vept\searrow0$, so that the ratio $\vept/\vepo^{\alpha_{d}}\searrow0$, for some suitable $\alpha_{d}$, then we expect the smoothing spatial operator to dominate the noise term and the solution to \eqref{lksop} or \eqref{heat} to tend in the limit, in a suitable norm, to the solution of the deterministic version of \eqref{lksop} or \eqref{heat} (obtained by setting $a\equiv0$);
\item (supercritical case) if, on the other hand, $\vepo\searrow0$ ``too fast'' compared to the speed at which $\vept\searrow0$, so that the ratio $\vept/\vepo^{\alpha_{d}}\nearrow\infty$, then we expect the  noise term to dominate the L-KS spatial operator.  The limiting behavior of our solution then bear some similarity to the more extreme---but simpler to see---limiting behavior when $\vept\nearrow\infty$ and $\vepo$ is kept fixed (thus $\vept/\vepo^{\alpha_{d}}\nearrow\infty$).  Assuming for simplicity that $\un\equiv0$ and $a$ is Lipschitz and bounded and bounded away from zero, we then have   
$$\E\lab U\rab^{2}=\vept^{2}\E\intrdzt K^{2}_{t-s;x,y}a^{2}(\usy)ds dy\to\infty$$
as $\vept\nearrow\infty$ for $d=1,2,3$, where $K$ is the fundamental solution of \eqref{lksop} or \eqref{heat} with $a\equiv0$.  We expect then that the supremum over space and time of the $L^{2}$ norm of the solution to grow to infinity in this supercritical case. 
\een
\fi

\subsubsection{\thmref{com}: from canonical L-KS SPDEs to nonlinear L-KS SPDEs via change of measure}  
At their core, our space-time change of measure theorems in \cite{Acom,Acom1,Acom2} are ``noise'' results that are independent of both the type and order of the SPDE under consideration.  This makes them conveniently adaptable to different SPDEs settings.  We use this fact to adapt our earlier change of measure results, from the second order equations in \cite{Acom,Acom1,Acom2} to the fourth order equations of this article, to transfer results and properties from the zero drift L-KS SPDE \eqref{nlkszd} (linear PDE part) to the nonzero-drift case \eqref{nlks} (nonlinear PDE part).  In addition, we use the same almost sure $L^2$ condition
on the drift/diffusion ratio as in our work \cite{Acom1,Acom2} to transfer uniqueness in law and establish law equivalence between solutions to \eqref{nlkszd} and \eqref{nlks}.  As observed in \cite{Acom1}, this is a much weaker condition than the traditional Novikov condition for change of measure; and this allows us to transfer results and properties from the canonical L-KS SPDEs \eqref{nlkszd} to many nonlinear L-KS SPDEs \eqref{nlks}, including the Swift-Hohenberg SPDE, driven by space-time white noise on subsets of $\Rp\times\Rd$, $d=1,2,3$. 

Now, we turn to the setting of our final main result of this paper.  Recall that we denote the zero-drift L-KS SPDE \eqref{nlkszd} by $\eLKSazun$ while the SPDE \eqref{nlks} is denoted by $\eLKSabun$.    We fix $T,L_{1},L_{2},L_{3}>0$, let $\T=[0,T]$, and we consider both equations $\eLKSazun$ and $\eLKSabun$ on the time-space domain $\T\times\S$, where either $\S=\Rd$ or $\S=\prod_{i=1}^{d}[0,L_{i}]$, $d=1,2,3$.  In the case $\S=\prod_{i=1}^{d}[0,L_{i}]$ the equations are supplemented with suitable boundary conditions\footnote{E.g., boundary conditions of Neumann type $\pa U/\pa n=\pa\D U/\pa n=0$ or Dirichlet type conditions  $U=\D U=0$ on $\pa\S$ and $d=1,2,3$.}, the nature of which is irrelevant to our change of measure results.  Also irrelevant to our results below is whether solutions are defined as mild kernel solutions like in \defnref{kskersoldef}---with appropriate modifications to $\KKSepthtx$ to account for the boundary conditions\footnote{E.g., in the Neumann (Dirichlet) case, the propagator ${ \e^{-|x-y|^2/2\i s}}/{{\lpa2\pi \i s \rpa}^{d/2}}$ in the definition of the $(\vep,\vth)$ L-KS kernel $\KKSepthtx$ \eqref{vepvthLKS} is replaced with the propagator with reflection (absorption) at $\pa\S$, respectively. } in the case $\S=\prod_{i=1}^{d}[0,L_{i}]$---or whether solutions are defined as weak; i.e., given in the test functions formulation (TFF).  For concreteness, and to also give the independently useful TFF for $\eLKSabun$, we take the TFF as our definition of solutions to $\eLKSabun$ for \thmref{com} and we now proceed to define it (see \remref{solequiv} below about the equivalence of the two formulations).  Let the Dirichlet test functions space be given by\footnote{Of course, the Dirichlet choice, which is assumed throughout the article whenever $\S=\prod_{i=1}^{d}[0,L_{i}]$, is without loss of generality and for concreteness only.  The Neumann (and other) boundary conditions are just as easily handled.}
\beq\lbl{tfspace}
\bsp
\Phi^{\infty}_{c,\mbox\tiny{Dir}}(\S;\R):=\bc
\lbr\varphi\in\C^{\infty}(\Rd;\R);\varphi=\D\varphi=0\mbox{ on }\pa\S\rbr;&\S=\prod\limits_{i=1}^{d}[0,L_{i}], \\
\C_{c}^{\infty}(\S;\R);&\S=\Rd,
\ec  
\end{split}
\eeq
where $d=1,2,3$.  
\begin{defn}[Test function solutions to $\eLKSabun$]\label{tfflksdef}
We say that the pair $(U,\sW)$ defined on the usual probability space
$\OFFtP$ is a test function solution to $\eLKSabun$ on $\Rp\times\S$ if $\sW$ is a space-time 
white noise on $\Rp\times\S$; the random field $U$ is progressively measurable, and with $U(0,x)=\unx$; and the pair $(U,\sW)$ satisfies the test function formulation:
\beq\lbl{tfflks}
\bsp
&\lpa U(t)-\un,\vph\rpa
=\int_{0}^{t}\lbk-\lpa U(s),\tfrac\vep8\lpa\lap+2\vth\rpa^{2}\vph\rpa +\lpa b(U(s)),\vph\rpa\rbk ds
\\&+ \int_{\S}\int_{0}^{t} a(\usy)\vph(y)\sW(ds\times dy); \ \forall \vph\in \Phi^{\infty}_{c,\mbox\tiny{Dir}}(\S;\R), t>0, \mbox{ a.s. } \P,
\end{split}
\eeq
where $\lpa\cdot,\cdot\rpa$ denotes the usual inner product on $L^{2}(\S;\R)$.  The test function solution is continuous if $U$ has continuous paths on $\Rp\times\S$.
Weak and strong---in the probability sense---solutions and uniqueness in law and pathwise uniqueness are defined in the usual way as in \defnref{kskersoldef}.
\end{defn}
\begin{rem}\lbl{solequiv}
We often simply say that $U$ is a test function solution to $\eLKSabun$ (weakly or strongly) to mean the
same thing as above.   As in Walsh's treatment of second order SPDEs (the top of p.~314 in \cite{W} and the discussion bedore it), it is straightforward to show the equivalence of the two formulations: kernel formulation in \eqref{ibtbapsol} (with spatial set $\S=\Rd$ or $\S=\prod_{i=1}^{d}[0,L_{i}]$) and test function formulation in \eqref{tfflks} under local boundedness assumptions on $a$ and $b$.
\end{rem}
We use $\lambda$ to denote the Lebesgue measure on $\sB\lpa\T\times\Rd\rpa$.  Also, for any function $u:\T\times\S$, we use the following notation for the drift/diffusion ratio function:
\beq\lbl{drdfratio} 
R_{u}(t,x):=\df{b(u(t,x))}{a(u(t,x))}; (t,x)\in\T\times\S.
\eeq
\bfr
\bthm[From canonical L-KS to nonlinear L-KS SPDEs on subsets of $\lbr\Rp\times\Rd\rbr_{d=1}^{3}$ via change of measure]\lbl{com}
Assume that either $\S=\Rd$ or $\S=\prod_{i=1}^{d}[0,L_{i}]$, $d=1,2,3$.  Suppose that the ratios $\RU$ and $\RV$ are in $L^2(\T\times\S,\lambda)$, almost surely,
whenever the continuous random fields $U$ and $V$ solve $($weakly or strongly\/$)$
$\eLKSazun$ and $\eLKSabun$, respectively, on $\T\times\S$.  Then,
\ben\rencomrom
\item uniqueness in law holds for $\eLKSabun$ iff uniqueness in law holds for $\eLKSazun;$ and
\item if uniqueness in law holds for $\eLKSazun$, $U$ is a  solution to $\eLKSazun$, and $V$ is a solution to $\eLKSabun$ on $\T\times\S;$ then the laws of $U$ and $V$ on $\sB\lpa\C(\T\times\S;\R)\rpa$ are equivalent (mutually absolutely continuous).
\een
In particular, let $\S=\prod_{i=1}^{d}[0,L_{i}]$, $d=1,2,3;$ and assume that $a(u)=\kappa\in\R\setminus\{0\}$,  $b(u)=\sum_{k=0}^{N}c_{k}u^{k}$ for $c_{k}\in\R$, $k=0,\ldots,N$, and $N\ge0$, and $u_0\in \mathrm{C}_{c}^{2,\gamma}(\S;\R)\mbox{ and nonrandom}$.  Then, the conclusions in $(i)$ and $(ii)$ above hold, uniqueness in law holds for $\eLKSabun$, and if $U$ and $V$ are continuous solutions to $\eLKSazun$ and $\eLKSabun$, respectively, on $\T\times\S$
\beqs\lbl{reglkssh}
\bsp
&V\in\H^{{\tf{4-d}{8}}^{-},{\lpa\tf{4-d}{2}\wedge 1\rpa}^{-}}\lpa\T\times\S;\R\rpa \mbox{ for $d=1,2,3$ a.s. }\\ \iff
&U\in\H^{{\tf{4-d}{8}}^{-},{\lpa\tf{4-d}{2}\wedge 1\rpa}^{-}}\lpa\T\times\S;\R\rpa \mbox{ for $d=1,2,3$  a.s.}
\end{split}
\eeqs
\ethm
\efr
The change of measure equivalence of \thmref{com}---between the canonical L-KS SPDE $\eLKSazun$ and a large class of nonlinear L-KS SPDEs $\eLKSabun$---immediately leads to uniqueness in law for the important special case ($N=2p-1$ $p\in\N$ and  $c_{2p-1}<0$) of the Swift-Hohenberg SPDE and its generalization: 
\beq\lbl{gSHcom}
\eLKSabun, \mbox{ with } b(u)=\sum_{k=0}^{2p-1}c_{k}u^{k}\mbox{ and with $p\in\N$ and  $c_{2p-1}<0$.}
\eeq
 It also says that the Swift-Hohenberg SPDE shares the same dimension-dependent local H\"older regularity with the canonical L-KS SPDE $\eLKSazun$, as in \thmref{com}.  Even more, the law of the solution of  the canonical L-KS SPDE is equivalent to that of the solution to the SH SPDE on $\sB\lpa\C(\T\times\S;\R)\rpa$, where $\S=\prod_{i=1}^{d}[0,L_{i}]$, $d=1,2,3.$      
\bfr
\bcor[Swift-Hohenberg uniqueness and law equivalence to the canonical L-KS SPDE]\lbl{SHcom}
Let $\S=\prod_{i=1}^{d}[0,L_{i}]$, $d=1,2,3,$  and assume that $a(u)=\kappa\in\R\setminus\{0\}$ and $u_0\in \mathrm{C}_{c}^{2,\gamma}(\S;\R)$ and nonrandom. The generalized Swift-Hohenberg SPDE \eqref{gSHcom} admits uniqueness in law and is law equivalent to $\eLKSazun$ on $\sB\lpa\C(\T\times\S;\R)\rpa;$ consequently,  it has the same H\"older regularity as the canonical L-KS SPDE $\eLKSazun$ on $\T\times\S$.\ecor
\efr 
\brm
We note that the conclusions of the last part of \thmref{com} hold also in the multiplicative noise case $a(u)=\kappa u$ and $b(u)=\sum_{k=1}^{N}c_{k}u^{k}$, where $\kappa\in\R\setminus\{0\}$ and $c_{i}\neq0$ for at least one $i\in\N$ (which covers the standard Allen-Cahn nonlinearity $u(1-u^{2})$ encountered in the SH equation).  We note here that all is needed is (1) uniqueness in law for $\eLKSazun$, which holds since the stronger pathwise uniqueness holds because $a(u)=\kappa u$ satisfies \eqref{lcnd} in \thmref{Holreg} and (2) the ratios $R_{U}$ are $R_{V}$ are clearly in $L^{2}(\T\times\S,\lambda)$ by the continuity assumption on $U$ and $V$ and the nonzero assumption on the constants $\kappa$ and the $c_{i}$'s\footnote{Of course we take $R_{U}{_{|_{U=0}}}:=\lim_{U\to0}R_{U}=\lim_{V\to0}R_{V}:=R_{V}{_{|_{V=0}}}=c_{1}/\kappa\mbox{ or }0$.}.
\erm
\subsection*{Acknowledgement}  I would like to sincerely thank the anonymous referee for his/her constructive comments which helped improve and clarify the presentation of this paper.

\section{A Harmonic connection between the L-KS and the BTBM kernels}\lbl{Foursec}
\subsection{Fourier transforms and $(\vep,\vth)$ L-KS PDEs links}
We start by obtaining the spatial Fourier transforms\footnote{We use the symmetric definition of the Fourier transform.  From a Physics point of view, the Fourier transform is taken over position to get energy.} for the Brownian-time Brownian motion (BTBM) and the $(\vep,\vth)$ L-KS kernels.  This reveals and captures both similarities and differences between both kernels and the PDEs corresponding to them.
\blm[Spatial Fourier transforms  of the BTBM and the $(\vep,\vth)$ L-KS kernels]\lbl{KSFT}
Let $\KBtx$ and $\KKSepthtx$ be the BTBM and $(\vep,\vth)$ L-KS kernels, respectively.
\ben\renrom
\item The spatial Fourier transform of $\KBtx$ is given by
\beq\lbl{BTBMF}
\FKBtxi=\lpa2\pi\rpa^{-\frac d2}\e^{\frac t8\lab\xi\rab^{4}}\lbk\frac{2}{\sqrt{\pi}}\int_{\frac{\sqrt{2t}\lab\xi\rab^{2}}{4}}^{\infty}\e^{-\tau^{2}}d\tau\rbk.
\eeq

\item The spatial Fourier transform of  $\KKSepthtx$ is given by
\beq\lbl{LKSF}
\FKKSepthtxi=\lpa2\pi\rpa^{-\frac d2}\e^{-\frac{\vep t}{8}\left( -2\vth+\lab\xi\rab^{2} \right) ^{2}};\ \vep>0,\ \vth\in\R.
\eeq
\een
\elm
\bpf
Starting with the BTBM kernel Fourier transform, we have
\beq\lbl{BTBMFcomp}
\bsp
\FKBtxi&=\lpa2\pi\rpa^{-\frac d2}\int_{\Rd}\lbk2\int_{0}^{\infty}\psxz\ptsz ds\rbk \e^{-\i\xi\cdot x}dx
\\&=\lpa2\pi\rpa^{-\frac d2}2\int_{0}^{\infty}\frac{\e^{-s^{2}/2t}}{\sqrt{2\pi t}}\e^{-\frac s2 \lab\xi\rab^{2}} ds
\\&=\lpa2\pi\rpa^{-\frac d2}\lbk\frac{2\e^{\frac t8\lab\xi\rab^{4}}}{\sqrt{\pi}}\int_{\frac{\sqrt{2t}\lab\xi\rab^{2}}{4}}^{\infty}\e^{-\tau^{2}}d\tau\rbk,
\end{split}
\eeq
proving part (i).  The Fourier transform of the $(\vep,\vth)$ L-KS kernel is now given by
\beq\lbl{lksFcomp}
\bsp
\FKKSepthtxi&=\lpa2\pi\rpa^{-\frac d2}\int_{\Rd}\lbk\int_{\R\setminus\{0\}}\df{\e^{i\vth s} \e^{-|x|^2/2\i s}}{{\lpa2\pi \i s \rpa}^{d/2}}\peptsz ds\rbk\e^{-\i\xi\cdot x}dx
\\&=\lpa2\pi\rpa^{-\frac d2}\int_{-\infty}^0\e^{-\frac{\i s}{2}\left( -2\vth+ \lab\xi\rab^{2} \right) }\peptsz ds+\int_{0}^\infty\e^{-\frac{\i s}{2}\left( -2\vth+ \lab\xi\rab^{2} \right) }\peptsz ds
\\&=\lpa2\pi\rpa^{-\frac d2}\int_{0}^\infty\frac{\e^{-s^{2}/2\vep t}}{\sqrt{2\pi\vep t}}\lbk\e^{-\frac{\i s}{2}\left( -2\vth+ \lab\xi\rab^{2} \right)}+\e^{\frac{\i s}{2}\left( -2\vth+ \lab\xi\rab^{2} \right)}\rbk ds
\\&=\lpa2\pi\rpa^{-\frac d2}\e^{-\frac{\vep t}{8}\left( -2\vth+\lab\xi\rab^{2} \right) ^{2}},
\end{split}
\eeq
completing the proof of \lemref{KSFT}.
\epf
\brm\lbl{kernelsrem}  The extra factor $\frac{2}{\sqrt{\pi}}\int_{\frac{\sqrt{2t}\lab\xi\rab^{2}}{4}}^{\infty}\e^{-\tau^{2}}d\tau$ in the BTBM transform \eqref{BTBMF} capture the memoryful property of the PDE \eqref{btppdedet} (the inclusion of $\un$) and the plus sign of the term $t|\xi|^{4}/8$ corresponds to that of the biLaplacian in \eqref{btppdedet}.     
\erm
Inverting the Fourier transform in \lemref{KSFT} we immediately get the more-convenient form for $\KKSepthtx$ in \eqref{iFTintro}, which can easily be verified to be a solution to the $(\vep,\vth)$ L-KS PDE in \eqref{vepvthlkspde}  with Dirac initial condition $\delta(x)$.  In particular, the special case $(\vep,\vth)=(1,1)$ confirms that our L-KS kernel $\KKStx$ in \eqref{ibtbapkernel} is the fundamental solution of the L-KS PDE in \eqref{lkspde}.  Let $u$ be given by 
\beq\lbl{detepthlksint}
u(t,x)=\int_{\Rd}\KKSepthtx\uny dy
\eeq 
and assume $\un$ satisfies the regularity conditions in \eqref{lcnd} (c).  The dominated convergence theorem plus a bit of analysis\footnote{See for example Lemma 2.1 in \cite{Aks}.  We leave the very similar details to the interested reader.} then give us that 
\beq\lbl{lkspdesecpf}
\bsp
\pat u(t,x)&=\int_{\Rd}\pat\KKSepthtx\uny dy=\int_{\Rd}-\tfrac\vep8\lpa\lap_{x}+2\vth\rpa^{2}\KKSepthtx\uny dy\\
&=-\tfrac\vep8\lpa\lap+2\vth\rpa^{2}\int_{\Rd}\KKSepthtx\uny dy=-\tfrac\vep8\lpa\lap+2\vth\rpa^{2}u(t,x)
\end{split}
\eeq
and $u(0,x)=\unx$.    Thus, we obtain the following theorem summarizing the PDEs connections.
\bthm\lbl{pdeconthm}
The $(\vep,\vth)$ L-KS kernel solves the initial value $(\vep,\vth)$ L-KS PDE 
\beq\lbl{vepvthlkspdedirac} 
\bsp
\frac{\partial u}{\partial t}&=-\tfrac\vep8\lpa\lap+2\vth\rpa^{2}u, t>0,x\in\Rd;\\
u(0,x)&=\delta(x), x\in\Rd.
\end{split}
\eeq
Moreover, if $u$ is given by \eqref{detepthlksint}, and $\un$ satisfies the condition in \eqref{lcnd} (c), then $u$ solves the $(\vep,\vth)$ L-KS PDE in \eqref{vepvthlkspdedirac} with $u(0,x)=\un(x)$.
\ethm
Setting $\vep=\vth=1$ in \eqref{lkspdesecpf} in the argument leading to \thmref{pdeconthm}, gives us an alternative proof of our Theorem 1.1 of \cite{Aks} connecting the linearized KS PDE \eqref{lkspde} to the L-KS kernel $\KKStx$.  On the other hand, setting $\vep=1$ and $\vth=0$ in \eqref{lkspdesecpf}; we get that the simpler kernel 
\beq\lbl{sfoker}
\Ksfotx:=\KKSepthoztx=\int_{-\infty}^0\df{e^{-|x-y|^2/2\i s}}{{\lpa2\pi \i s \rpa}^{d/2}}\ptzs ds+\int_{0}^\infty\df{e^{-|x-y|^2/2\i s}}{{\lpa2\pi \i s \rpa}^{d/2}}\ptzs ds,
\eeq
obtained by removing the angle $\e^{\i s}$ from the L-KS kernel $\KKStx$ in \eqref{ibtbapkernel}, is the fundamental solution of the simpler fourth order PDE 
\begin{equation} \label{sfopde}
 \begin{cases} \displaystyle\frac{\partial u}{\partial t}=
-\frac18\lap^2u, & (t,x)\in(0,+\infty )\times\Rd;
\cr u(0,x)=\delta(x), & x\in\Rd
\end{cases}
\end{equation}
as was shown for the case $d=1$ in Hochberg and Orsinger \cite{HO} (see also the  different approach in Funaki \cite{Fun}, also for $d=1$).  Clearly, the Fourier transforms $\FKKStxi$ and $\FKsfotxi$ of $\KKStx$ and $\Ksfotx$, and their inverses are now given as an immediate corollary to \lemref{KSFT}.  Taking $(\vep,\vth)=(1,1)$ and $(\vep,\vth)=(1,0)$, respectively in \lemref{KSFT} (ii) and using a dominated convergence argument, we get
\bcor\lbl{FTandinvFTcor}
\beq\lbl{FTsfoker}
\bsp
\FKKStxi=\df{\e^{-\frac{t}{8}\left( -2+\lab\xi\rab^{2} \right) ^{2}}}{\lpa2\pi\rpa^{\frac d2}},\ 
\KKStx&=(2\pi)^{-d}\int_{\Rd}\e^{-\frac{t}{8}\left( -2+\lab\xi\rab^{2} \right) ^{2}}\e^{\i\xi\cdot x} d\xi \\ 
&=(2\pi)^{-d}\int_{\Rd}\e^{-\frac{t}{8}\left( -2+\lab\xi\rab^{2} \right) ^{2}}\cos\lpa{\xi\cdot x}\rpa d\xi;
\\
\FKsfotxi=\frac{\e^{-\frac{t}{8}\lab\xi\rab^{4}}}{\lpa2\pi\rpa^{\frac d2}},\ 
\Ksfotx&=(2\pi)^{-d}\int_{\Rd}\e^{-\frac{t}{8}\lab\xi\rab^{4}}\e^{\i\xi\cdot x} d\xi\\
&=(2\pi)^{-d}\int_{\Rd}\e^{-\frac{t}{8}\lab\xi\rab^{4}}\cos\lpa{\xi\cdot x}\rpa d\xi.
\end{split}
\eeq
\ecor

\subsection{A revealing kernels $L^{2}$ energy}
To understand why the L-KS and the BTBM kernels $\KKStx$ and $\KBtx$ have very similar regularizing effects on the L-KS SPDE \eqref{nlkszd} above (with $(\vep,\vth)=(1,1)$) and the BTBM SIE introduced in \cite{Abtbmsie} (and obtained from \eqref{ibtbapsol} by replacing $\KKSepthtx$ with $\KBtx$ and setting $b\equiv0$), we first observe that the regularity of the L-KS PDE \eqref{lkspde} is dictated by the bi-Laplacian term and that the family $$\lbr\KKSepthtx\rbr_{\vep>0,\vth\in\R}$$ of all $(\vep,\vth)$ L-KS kernels in \eqref{vepvthLKS} and \eqref{iFTintro}---including $\KKStx$ and $\Ksfotx$---share the same regularizing effect on the L-KS SPDE \eqref{nlkszd}.  

As we will see shortly, the $L^{2}$ quantity 
\beq\lbl{L2ofsfoker}
 \intrd\lab\Ksfotx\rab^2dx=\intrd\lab\FKsfotxi\rab^2d\xi=\intrd\frac{\e^{-\frac{t}{4}\lab\xi\rab^{4}}}{\lpa2\pi\rpa^{d}}d\xi=C_{d}t^{-d/4}; d=1,2,3,
 \eeq
 for the $(\vep,\vth)=(1,0)$ L-KS kernel $\Ksfotx$---where we used the Parseval-Plancherel theorem and where $C_{d}$ is a dimension dependentt constant\footnote{$C_{1}={1}/{2\Gamma\lpa\frac34\rpa}$, $C_{2}=1/4\sqrt{\pi}$, and $C_{3}=\Gamma\lpa\frac34\rpa/\pi^{2}\sqrt8$.}---is key to understanding the regularity of our L-KS SPDE \eqref{nlkszd}.   By the above discussion (see also \lemref{L2} below), it is clear that $\intrd\lab\KKStx\rab^2dx$ is of the same order\footnote{In fact, in $d=2$ $$\intrd\lab\FKKStxi\rab^2d\xi=\lbk1+\psi(\sqrt{t})\rbk\intrd\lab\FKsfotxi\rab^2d\xi,\ t>0$$
 where $\psi(u):=(2/\sqrt{\pi})\int_{0}^{u}\e^{-r^{2}}dr$.  See also \lemref{L2} below.}.   On the other hand as was shown in Lemma 2.2 in \cite{Abtbmsie}, there is a dimension dependent constant $c_{d}$ such that 
 \beq\lbl{L2btbmker}
\intrd\lbk\KBtx\rbk^2dx=c_{d}t^{-d/4};\quad t>0,\ d=1,2,3.
 \eeq
 Equations \eqref{L2ofsfoker} and \eqref{L2btbmker} are the fundamental analytic reason why the regularity for our L-KS SPDE in our first result \thmref{Holreg} above is the same as that of the BTBM SIE in Theorem 1.1 of \cite{Abtbmsie}, albeit here we have real solutions to a negative bi-Laplacian equation and the BTBM SIE in \cite{Abtbmsie} has real solutions to a positive bi-Laplacian equation with memory (see \cite{Abtbmsie}).
\section{Proof of \thmref{Holreg}}
Since both $\vep>0$ and $\vth\in\R$ are fixed in \thmref{Holreg}, and since all the main conclusions are unaffected by the specific values of $\vep>0$ and $\vth\in\R$; we will simplify our notation and exposition by assuming  throughout this section (and its subsections)---without loss of generality---that either $(\vep,\vth)=(1,1)$ (capturing the general biLaplacian, Laplacian, and zero order term case) or $(\vep,\vth)=(1,0)$ (the biLpalcian term, without the lower order terms, case)\footnote{It should be clear that our methods extend with only minor notational changes to any fixed values for $\vep>0$ and $\vth\in\R$.  The case $(\vep,\vth)=(1,0)$ is the simplest representative case, and we include it explicitly in this subsection since it is useful in \lemref{L2} to obtain the fundamental $L^{2}$ estimates for the more interesting $(\vep,\vth)=(1,1)$ case.}.  
\subsection{Key regularity estimates for the L-KS kernel}\lbl{kerregest}
Here, we prove several $L^{2}$ estimates\footnote{\lemref{L2} is stated only for $0<t\le T$, since we only need it for intervals $[0,T]$.  In fact, for $d=2$, we show that the estimates hold, with the same constants $C_{l}^{2}$ and $C_{u}^{2}$, for all $t>0$.} on the L-KS kernel and its temporal and spatial differences that are key in proving our regularity results in \thmref{Holreg}.  Again, these fundamental estimates for the L-KS kernel are very similar to those for the BTBM density in the corresponding estimates in \cite{Abtbmsie}, but  the proofs proceed differently due to the oscillatory nature of the modified propagator part of the  L-KS kernel.
\begin{lem}[Kernel's $L^{2}$]\label{L2}  Fix any arbitrary  $T>0$.  There are constants $C_{l}^{(d)}$ and $C_{u}^{(d)}$ depending only on the spatial dimension $d$ and $T$ and a constant $C_{d}$ depending only on $d$ such that
\beq\lbl{L2wotint}
\bsp
&\intrd\lab\Ksfotx\rab^2dx=C_{d}t^{-d/4}; \mbox{ and }
\\ C_{l}^{(d)}t^{\tf{-d}{4}}\le&\intrd\lab\KKSsx\rab^2dx \le C_{u}^{(d)} t^{\tf{-d}{4}};
\end{split}
\eeq
for $0<t\le T,\ d\in\{1,2,3\}$ and hence
\beq\lbl{L2wtint}
\bsp
&\int_0^t\intrd\lab\Ksfosx\rab^2dx ds=C_{d}t^{\tf{4-d}{4}}; \mbox{ and }
\\ C_{l}^{(d)}t^{\tf{4-d}{4}}\le&\int_0^t\intrd\lab\KKSsx\rab^2dx ds\le C_{u}^{(d)} t^{\tf{4-d}{4}};
\end{split}
\eeq
for $0<t\le T,\ d\in\{1,2,3\}.$
\end{lem}
\bpf The equalities in \eqref{L2wotint} and in \eqref{L2wtint} follow immediately from \eqref{L2ofsfoker}.  Using the Parseval-Plancherel theorem, we have 
\beq\lbl{papl}
\bsp
\intrd\lab\KKSsx\rab^2dx=\intrd\lab\FKKSsxi\rab^2d\xi=\lpa2\pi\rpa^{-d}\intrd\e^{-\frac{s}{4}\left( -2+\lab\xi\rab^{2} \right) ^{2}}d\xi,
\end{split}
\eeq
for every $s>0$.  Let $d=2$ and $\psi(u):=(2/\sqrt{\pi})\int_{0}^{u}\e^{-r^{2}}dr$.  We then have
\beq\lbl{2d}
\bsp
\frac{1}{4\sqrt{\pi s}}&\le
\lpa2\pi\rpa^{-2}\intrs\e^{-\frac{s}{4}\left( -2+\lab\xi\rab^{2} \right) ^{2}}d\xi=\frac{1+\psi(\sqrt{s})}{4\sqrt{\pi}}\df{1}{\sqrt{s}}\\&=\lbk1+\psi(\sqrt{s})\rbk\intrd\lab\FKsfosxi\rab^2d\xi\le\frac{1}{2\sqrt{\pi s}}
\end{split}
\eeq 
and the assertions in \eqref{L2wotint} and its immediate consequence \eqref{L2wtint} are established for $d=2$.  

For dimensions $d=1,3$, we get the desired estimates by comparing $\intrd\lab\FKKStxi\rab^2d\xi$ with $\intrd\lab\FKsfotxi\rab^2d\xi$ (see \eqref{FTsfoker} and \eqref{L2ofsfoker} above).  To start, we use  \eqref{FTsfoker} and observe that
 \beq\lbl{asym}
 \lim_{t\searrow0} \frac{\ds\intrd\lab\FKKStxi\rab^2d\xi}{\ds\intrd\lab\FKsfotxi\rab^2d\xi}
 =\lim_{t\searrow0}\frac{\ds\intrd\e^{-\frac{t}{4}\left( -2+\lab\xi\rab^{2} \right) ^{2}}d\xi}{\ds\intrd\e^{-\frac{t}{4}\lab\xi\rab^{4}}d\xi}=1; \quad d=1,2,3.
 \eeq
From \eqref{FTsfoker}, \eqref{L2ofsfoker}, and \eqref{asym}, we then easily have 
\beq\lbl{supinfconst}
\bsp
&C_{\mbox{\tiny min}}^{(d)}:=\inf_{0<t\le T}\frac{\ds\int_{\Rd}\lab\FKKStxi\rab^2d\xi}{\ds\int_{\Rd}\lab\FKsfotxi\rab^2d\xi}=\inf_{0<t\le T}\frac{\ds\intrd\e^{-\frac{t}{4}\left( -2+\lab\xi\rab^{2} \right) ^{2}}d\xi}{\ds\intrd\e^{-\frac{t}{4}\lab\xi\rab^{4}}d\xi}>0\mbox{ and }\\
&C_{\mbox{\tiny max}}^{(d)}:=\sup_{0<t\le T}\frac{\ds\int_{\Rd}\lab\FKKStxi\rab^2d\xi}{\ds\int_{\Rd}\lab\FKsfotxi\rab^2d\xi}=\sup_{0<t\le T}\frac{\ds\intrd\e^{-\frac{t}{4}\left( -2+\lab\xi\rab^{2} \right) ^{2}}d\xi}{\ds\intrd\e^{-\frac{t}{4}\lab\xi\rab^{4}}d\xi}<\infty,
\end{split}
\eeq 
for $d-1,2,3$.  So, for $d=1,3$, and $0<s\le T$, we use the Parseval-Plancherel theorem together with \eqref{supinfconst} and \eqref{L2ofsfoker} to get the desired lower and upper bounds as follows: 
\beq\lbl{L2wotint13dlub}
\bsp
\int_{\Rd}\lab\KKSsx\rab^2dx&=\int_{\Rd}\lab\FKKSsxi\rab^2d\xi \ge C_{\mbox{\tiny min}}^{(d)}\int_{\Rd}\lab\FKsfosxi\rab^2d\xi= C_{\mbox{\tiny min}}^{(d)}C_{d} s^{-\tf{d}{4}},\\
\int_{\Rd}\lab\KKSsx\rab^2dx&=\int_{\Rd}\lab\FKKStxi\rab^2d\xi \le C_{\mbox{\tiny max}}^{(d)}\int_{\R}\lab\FKsfotxi\rab^2d\xi= C_{\mbox{\tiny max}}^{(d)}C_{d} s^{-\tf{d}{4}}.
\end{split}
\eeq
The assertions in \eqref{L2wotint} and its immediate consequence \eqref{L2wtint} are thus established for $d=1,3$ and rhe proof is complete.
\epf
\brm\lbl{lksorsl2relation}
In $d=1$, there is a critical $t_{c}>1$ such that\footnote{$t_{c}\approx1.506188$.  It is interesting to note that this is only a one-dimensional phenomenon (see \eqref{2d} and \eqref{3dim}).} 
\beq\lbl{1dim}
\bc\ds\int_{\R}\lab\FKsfotxi\rab^2d\xi<\int_{\R}\lab\FKKStxi\rab^2d\xi,&t<t_{c}\vspace{1mm}\cr
\ds\int_{\R}\lab\FKsfotxi\rab^2d\xi\ge\int_{\R}\lab\FKKStxi\rab^2d\xi,&t\ge t_{c},
\ec
\eeq
with equality at $t=t_{c}$.  If $T\le t_{c}$, then using \eqref{L2ofsfoker} and \eqref{1dim} the lower bound of \eqref{L2wotint} immediately holds with $C_{l}^{(1)}=C_{1}=1/{2\Gamma\lpa\frac34\rpa}$, where $C_{1}$ is the constant in \eqref{L2ofsfoker} for $d=1$.   On the other hand, as in the case $d=2$ (see \eqref{2d} above), when $d=3$ we have 
\beq\lbl{3dim}
\int_{\R^{3}}\lab\FKsfotxi\rab^2d\xi<\int_{\R^{3}}\lab\FKKStxi\rab^2d\xi; \ t>0,
\eeq
which, when combined with \eqref{L2ofsfoker}, gives us the lower bound with the constant $C_{l}^{(3)}=C_{3}=\Gamma\lpa\frac34\rpa/\pi^{2}\sqrt8$.  
\erm
\begin{lem}[Kernel's $L^{2}$ temporal difference]\label{temp}  Fix any arbitrary  $T>0$.
There are constants $\tilde{C}_{u}^{(d)}$, depending only on $d$ and $T$ such that 
\beq\lbl{tmpkernelup}
\int_0^t\intrd{\lab\KKStsx - \KKSrsx\rab}^2 dx ds \leq \tilde{C}_{u}^{(d)}(t-r)^{\tf{4-d}{4}};0<r<t\le T, d=1,2,3,
\eeq
with the convention that $\KKStx=0$ if $t<0$.  The same estimate holds, with possibly different constants, when replacing $\KKStx$ with $\Ksfotx$.
\end{lem}
\bpf
Throughout the proof, unless otherwise specified, the spatial dimension $d\in\lbr1,2,3\rbr$. For $u,v>0$ let 
\beq\lbl{Ktildedef}
\tilde{\K}_{u+v}^{(d)}=(2\pi)^{-d}\intrd\e^{-\frac{u}{8}\left( -2+\lab\xi\rab^{2} \right) ^{2}}\e^{-\frac{v}{8}\left( -2+\lab\xi\rab^{2} \right) ^{2}}d\xi.
\eeq
By the Parseval-Plancherel theorem, we have 
\beq\lbl{tmpkernelPP}
\bsp
&\int_0^t\intrd{\lab\KKSsptmrx - \KKSsx\rab}^2 dx ds =
\int_0^t\intrd{\lab\FKKSsptmrxi - \FKKSsxi\rab}^2 d\xi ds
\\&=\int_0^t\lbk\tilde{\K}_{2[s+(t-r)]}^{(d)}-2\tilde{\K}_{2s+(t-r)}^{(d)}+\tilde{\K}_{2s}^{(d)} \rbk ds
\\&=\lbk\int_{0}^{\tf{t-r}{2}}\tilde{\K}_{2s}^{(d)} ds-\int_{\tf{t-r}{2}}^{t-r}\tilde{\K}_{2s}^{(d)} ds-\int_{t}^{t+\tf{t-r}{2}}\tilde{\K}_{2s}^{(d)} ds
+\int_{t+\tf{t-r}{2}}^{2t-r}\tilde{\K}_{2s}^{(d)} ds\rbk.
\end{split}
\eeq
It is clear from \eqref{Ktildedef} that $\tilde{\K}_{2s}^{(d)}$ is decreasing in $s$.  Thus, the sum of the last three terms of \eqref{tmpkernelPP} is $\le0$ and we have
\beq\lbl{tmpkernelub}
\bsp
&\int_0^t\intrd{\lab\KKSsptmrx - \KKSsx\rab}^2 dx ds \le \int_{0}^{\tf{t-r}{2}}\tilde{\K}_{2s} ^{(d)} ds
\\&=\int_0^{\tf{t-r}{2}}\intrd\lab\KKSsx\rab^2dx ds\le \tilde{C}_{u}^{(d)}(t-r)^{\tf{4-d}{4}};0<r<t\le T,d\in\{1,2,3\},
\end{split}
\eeq
 where we used the definition of $\tilde{\K}_{2s} ^{(d)}$ in \eqref{Ktildedef}, Parseval-Plancherel theorem, and \lemref{L2}.  The proof of the simpler $\Ksfotx$ case follows the same steps, with obvious trivial changes, and will be omitted.  The lemma is established\footnote{The constants $\tilde{C}_{u}^{(d)}=\lbk {2^{\tf{d-4}{4}}}\rbk{C}_{u}^{(d)}$, where the constants ${C}_{u}^{(d)}$ are those in \lemref{L2}.}.  
 \epf

\begin{lem}[Kernel's $L^{2}$ spatial difference]\label{space}
For $d\in\{1,2,3\}$, there are intervals $I_{1}=(0,1]$, $I_{2}=(0,1)$, and $I_{3}=(0,1/2);$
positive numbers $\lbr\alpha_{d}\in I_{d}\rbr_{d=1}^{3}$; constants $\lbr C^{(d)}_{u}\rbr_{d=1}^{3}$ depending only on $d$ and $\alpha_{d}\in I_{d}$ such that 
\beq\lbl{sptlkernel}
\ds\int_0^t\intrd {\lab\KKSsx - \KKSsxpz\rab}^2 dx ds \leq C^{(d)}_{u} |z|^{2\alpha_{d}}(1\vee\tfrac t4); \forall\alpha_{d}\in I_{d}, t>0,
\eeq
where $0<C^{(d)}_{u}<\infty$ for every $\alpha_{d}\in I_{d}$ for $d=1,2,3$.  The same estimate holds, with possibly different constants, when replacing $\KKStx$ with $\Ksfotx$. 
\end{lem}
\bpf   We first observe from \eqref{lksFcomp} that 
\beq\lbl{shiftfour}
\FKKSsxipz=\lpa2\pi\rpa^{-\frac d2}\e^{-\frac{t}{8}\left( -2+\lab\xi\rab^{2} \right) ^{2}}\e^{\i z\cdot \xi}.
\eeq
Suppose $d\in\{1,2,3\}$, and let $\B^{d}_{\sqrt{2}}:=\lbr\xi\in\Rd;\lab\xi\rab<\sqrt2\rbr$.  Again, the Parseval-Plancherel theorem tells us that the quantity we want to estimate is 
\beq\lbl{spatdiff}
\bsp
&\int_0^t\intrd {\lab\FKKSsxi - \FKKSsxipz\rab}^2 d\xi ds
\\&=(2\pi)^{-d}\int_0^t\intrd\lab\e^{-\frac{s}{8}\left( -2+\lab\xi\rab^{2} \right) ^{2}}\lbk1-\e^{\i z\cdot \xi}\rbk\rab^{2}d\xi ds
\\&=2(2\pi)^{-d}\int_0^t\intrd \e^{-\frac{s}{4}\left( -2+\lab\xi\rab^{2} \right) ^{2}}\lbk1-\cos\lpa z\cdot\xi\rpa\rbk d\xi ds
\\&=8(2\pi)^{-d}\int_{\B^{d}_{\sqrt{2}}} \lbk\tfrac{1-\e^{-\frac{t}{4}\left( -2+\lab\xi\rab^{2} \right) ^{2}}}{\left( -2+\lab\xi\rab^{2} \right) ^{2}}\rbk\lbk1-\cos\lpa z\cdot\xi\rpa\rbk d\xi
\\&+8(2\pi)^{-d}\int_{\Rd\setminus\lpa\B^{d}_{\sqrt{2}}\cup\partial\B^{d}_{\sqrt{2}}\rpa} \lbk\tfrac{1-\e^{-\frac{t}{4}\left( -2+\lab\xi\rab^{2} \right) ^{2}}}{\left( -2+\lab\xi\rab^{2} \right) ^{2}}\rbk\lbk1-\cos\lpa z\cdot\xi\rpa\rbk d\xi.
\end{split}
\eeq
We make use of the following two sets of elementary inequalities for all $d\ge1$, the first of  which uses the Cauchy-Schwarz inequality to obtain the last bound
\beq\lbl{elemineq}
\bsp
1-\cos\lpa z\cdot\xi\rpa&\le 2\lpa1\wedge\lab z\cdot\xi\rab^{2\alpha}\rpa\le2\lpa1\wedge\lab z\rab^{2\alpha}\lab\xi\rab^{2\alpha}\rpa;\ 0<\alpha\le1,
\\ \frac{1-\e^{-\frac{t}{4}\left(-2+\lab\xi\rab^{2} \right)^{2}}}{\lpa-2+\lab\xi\rab^{2}\rpa^{2}}&
\le(1\vee\tfrac t4)\wedge\frac{\lpa1\vee\tfrac t4\rpa\lpa1-\e^{-\lpa-2+\lab\xi\rab^{2}\rpa^{2}}\rpa}{\left( -2+\lab\xi\rab^{2}\rpa^{2}}; \  t\ge0.
\end{split}
\eeq

 We now treat the cases $d=1,2,3$ separately.
Using \eqref{spatdiff}, \eqref{elemineq}, and changing to polar coordinates in $d=2$ and to spherical coordinates in $d=3$ we can bound our desired quantity  
$$\int_0^t\intrd {\lab\FKKSsxi - \FKKSsxipz\rab}^2 d\xi ds$$ 
from above by
\beq\lbl{onedspat}
\bsp
&\frac{16(1\vee\tfrac t4)}{2\pi}|z|^{2\alpha}\lbk\int_{-\sqrt{2}}^{\sqrt{2}} \lab\xi\rab^{2\alpha} d\xi+2\int_{\sqrt{2}}^{\infty} \tfrac{1-\e^{-\left(-2+\lab\xi\rab^{2} \right)^{2}}}{\lpa-2+\lab\xi\rab^{2}\rpa^{2}}\lab\xi\rab^{2\alpha} d\xi\rbk
\\&\le C^{(1)} |z|^{2\alpha};\ 0<\alpha\le1 \mbox{ for }d=1,
\end{split}
\eeq
\beq\lbl{twodspat}
\bsp
&\frac{16(1\vee\tfrac t4)}{(2\pi)^{2}}|z|^{2\alpha}\lbk\int_{0}^{2\pi}\int_{0}^{\sqrt{2}} r^{2\alpha} rdr d\theta+\int_{0}^{2\pi}\int_{\sqrt{2}}^{\infty}
\tfrac{1-\e^{-\left(-2+r^{2}\right)^{2}}}{\lpa-2+r^{2}\rpa^{2}}r^{2\alpha} rdr d\theta\rbk
\\&\le C^{(2)} |z|^{2\alpha};\ 0<\alpha<1 \mbox{ for }d=2,
\end{split}
\eeq
and 
\beq\lbl{threedspat}
\bsp
&\frac{16(1\vee\tfrac t4)|z|^{2\alpha}}{(2\pi)^{3}}\lbk\int_{0}^{\pi}\int_{0}^{2\pi}\int_{0}^{\sqrt{2}} r^{2\alpha} r^{2}\sin(\vartheta)dr d\theta d\vartheta\right.
\\&\hspace{1in}+\left.\int_{0}^{\pi}\int_{0}^{2\pi}\int_{\sqrt{2}}^{\infty}
\tfrac{1-\e^{-\left(-2+r^{2}\right)^{2}}}{\lpa-2+r^{2}\rpa^{2}}r^{2\alpha} r^{2}\sin(\vartheta)dr d\theta d\vartheta\rbk
\\&\le C^{(3)} |z|^{2\alpha};\ 0<\alpha<\tfrac12 \mbox{ for }d=3.
\end{split}
\eeq
In particular, when $d=1$, $\alpha$ may be taken to be $1$; in $d=2$, $\alpha\in(0,1)$; and in $d=3$, $\alpha\in(0,1/2)$. 
Our dimension-dependent upper bound constant $C^{(d)}$ is independent of $t$ if $t\le4$ and increases with $t$ if $t>4$.  The proof of the simpler $\Ksfotx$ case follows the same steps, with obvious trivial changes, and will be omitted.
\epf
\subsection{Finishing the proof of \thmref{Holreg}}
We now complete the proof of \thmref{Holreg} by first establishing the H\"older regularity results without imposing any Lipschitz conditions on $a$, assuming the $L^{p}$ boundedness of solutions on $\T\times\Rd$.  We then add a Lipschitz condition on $a$ and obtain the strong (stochastically) existence and uniqueness result for the L-KS SPDE \eqref{nlkszd}, together with the $L^{p}$ boundedness assumed before; we thus obtain the H\"older regularity with no $L^{p}$ boundedness assumptions\footnote{As we mentioned in \remref{noliprem} the existence of lattice limit solutions along with the regularity results in \thmref{Holreg} (including both $L^{p}$ boundedness on $\T\times\Rd$ and H\"older regularity) can be proven under the weaker non-Lipschitz conditions \eqref{nlcnd}, as we did in \cite{Atfhosie,Abtbmsie}.}.  With \lemref{L2}--\lemref{space} in hand, the rest of the proof of  \thmref{Holreg} is a straightforward adaptation of our corresponding arguments in \cite{Abtbmsie} to our setting here.  For the convenience of the reader; we state the needed results  below, referring to \cite{Abtbmsie} for the details and noting the minor notational changes. Without loss of generality, assume for the remainder of the section that $(\vep,\vth)=(1,1)$ and that $U$ solves  the L-KS SPDE \eqref{nlkszd}.

Recalling that the initial data $\un$ is assumed deterministic and writing $U$ in the kernel formulation \eqref{ibtbapsol} in terms of  its deterministic and random parts  $\utx={U}_D(t,x)+{U}_R(t,x)$, we note that the deterministic part ${U}_D(t,x)=\intrd\KKStxy \uny dy$ is $\C^{1,4}(\Rp\times\Rd;\cC)$ smooth in time and space, under the assumptions on $\un$, since it is a classical solution to the LKS PDE \eqref{nlkszd} for $(\vep,\vth)=(1,1)$.  So, it suffices to get the H\"older regularity of the random part.  We now give estimates on the spatial and temporal differences of the random part ${U}_R$.  To get straight to these important regularity estimates, we first assume that 
\beq\lbl{momcnd}
M_q(t)=\sup_{x\in\Rd}\E|U(t,x)|^{2q}\le K_{T,q}<\infty; \quad t\in\T=[0,T],q\ge1.
\eeq
 \begin{lem}[Spatial and temporal differences $L^{2q}$ estimates] \lbl{sptdiff}
Assume that \eqref{lcnd} and \eqref{momcnd} are in force.  There exists a constant $\tilde{C}_{d}$ depending only on $q$, $\max_x |u_0(x)|$, the spatial dimension
$d\in\{1,2,3\}$, $\alpha_{d}$, and $T$ such that
\[
\mathbb{E}\left| {U}_R(t,x) - {U}_R(t,y\right|^{2q}
\le \tilde{C}_{d} |x-y|^{2q\alpha_{d}};\ \alpha_{d}\in I_{d},
\]
for all $x,y \in \Rd$, $t\in\T$, and $d\in\{1,2,3\};$ where $\alpha_{d}$ and $I_{d}$ are as in \lemref{space}.  

Also, there exists a constant $\bar{C}_{d}$ depending only on $q$, $\max_x|u_0(x)|$, the spatial dimension
$d\in\{1,2,3\}$, and $T$ such that
\[\mathbb{E}\left| {U}_R(t,x) - {U}_R(r,x) \right|^{2q}
\leq \bar{C}_{d}\lab t-r\rab^{\tf{(4-d)q}{4}},\] for all $x \in \Rd$, for all $t,r \in\mathbb{T}$, and for $1\le d\le 3$.
\end{lem}
\begin{pf} 
The proof follows the same exact steps as in the proofs of Lemma 2.5 and Lemma 2.6 in \cite{Abtbmsie}, replacing the kernel there with the L-KS kernel here and using \lemref{L2}--\lemref{space} here in place of the corresponding ones there.  
\end{pf}
We now have the desired H\"older regularity result as the following corollary.
\bcor[H\"older regularity]\lbl{Holconcl}
Assume that $(U,\sW)$ is an L-KS kernel solution to \eqref{nlkszd} on $\{\Rp\times\Rd\}_{d-1}^{3}$.  Suppose further that the $L^{p}$ boundedness in \eqref{momcnd} holds.  Then
$$U\in\H^{{\tf{4-d}{8}}^{-},{\lpa\tf{4-d}{2}\wedge 1\rpa}^{-}}\lpa\T\times\Rd;\cC\rpa;\mbox{ for  $d=1,2,3$.}$$
almost surely.
\ecor
\bpf
It suffices to prove it for the random part.  By \lemref{sptdiff} we easily have
\beq\lbl{directbtpsiehldr}
\bc
\E\lab{U}_R(t,x) - {U}_R(t,y)\rab^{2n+2d}\le C_{d}\lab x-y\rab^{(2n+2d)\alpha_{d}},\\
\E\lab{U}_R(t,x) - {U}_R(r,x)\rab^{2m+4d}\le \bar{C}_{d}\lab t-r\rab^{\tf{(4-d)(m+2d)}{4}},
\ec
\eeq
for $d=1,2,3$.  Thus, by standard results, the spatial H\"older exponent is
$\gamma_s\in\lpa 0,\tf{2(n+d)\alpha_{d}-d}{2n+2d}\rpa$ and the temporal exponent is $\gamma_t\in\lpa0,\tf{m\lpa1-d/4\rpa+d(1-d/2)}{2m+4d}\rpa$
$\forall m,n$.  Taking the limits as $m,n\to\infty$, we get
$\gamma_t\in\lpa0,\tf{4-d}{8}\rpa$ and $\gamma_s\in\lpa0,\alpha_{d}\rpa$, for $1\le d\le 3$ and $\alpha\in I_{d}$ as in \lemref{space}.  The proof is now complete. 
\epf
The final piece needed for the proof of \thmref{Holreg} is now given by the following Lemma, which also removes  the $L^{p}$ boundedness assumption \eqref{momcnd} by asserting that it automatically holds under the conditions \eqref{lcnd}.
\blm[Existence, uniqueness, and $L^{p}$ boundedness]\lbl{exunlp}  Under the conditions in \eqref{lcnd}, there exists a strong and pathwise unique solution to the L-KS SPDE \eqref{nlkszd} on $\Rp\times\Rd$ that is $L^{p}(\Omega)$-bounded on $\T\times\Rd$, for every $p\ge2$ and every $d=1,2,3$. 
\elm
\bpf
The proof follows exactly the same steps as the proof on pp. 27--29 in \cite{Abtbmsie} for the BTBM SIE, with now obvious and minor changes from the BTBM setting of \cite{Abtbmsie} to our L-KS setting here, we omit the details and point the interested reader to \cite{Abtbmsie} for the specifics.   The proof of \thmref{Holreg} is now complete.
\epf
\section{Proof of \thmref{intense}}  
We now turn to the proof of \thmref{intense}.  Again,  without loss of generality, it is enough to fix\footnote{Again, \thmref{intense} holds for any fixed value of $\vth\in\R$.  Fixing $\vth$ value to $1$ is for convenience and for simplifying notation and exposition in the proof.only.   The case $\vth=1$ captures the general case of the 4th order biLaplacian term together with the second and zero order terms.} $\vth=1$.  We first need the $\vepo$-time-scaled L-KS Kernel $\KKSeptx$, which---upon taking $\vth=1$ in \eqref{iFTintro}--reduces to
\beq\lbl{iFTformepeq}
\bsp
\KKSeptx=\KKSepotx&=(2\pi)^{-d}\int_{\Rd}\e^{-\frac{\vepo t}{8}\left( -2+\lab\xi\rab^{2} \right) ^{2}}\e^{\i\xi\cdot x} d\xi,
\\&=(2\pi)^{-d}\int_{\Rd}\e^{-\frac{\vepo t}{8}\left( -2+\lab\xi\rab^{2} \right) ^{2}}\cos\lpa{\xi\cdot x}\rpa d\xi
\end{split}
\eeq
and which, by \thmref{pdeconthm}, solves the L-KS PDE 
\begin{equation} \label{eplkspde}
 \begin{cases} {\ds\frac{\partial u}{\partial t}}=
-\frac\vepo8\lap^2u-\frac\vepo2\lap u-\frac\vepo2u, & (t,x)\in(0,+\infty )\times\Rd;
\cr u(0,x)=\delta(x); & x\in\Rd.
\end{cases}
\end{equation}
Also, exactly as in \lemref{L2} above, $\KKSeptx$ satisfies the bounds\footnote{Of course, the case $\vth=0$ satisfies \eqref{eplkspde}, with only the biLaplacian term, without the Laplacian and without the zero order terms; and satisfies \eqref{eplksl2} with equality.}.
\beq\lbl{eplksl2}
\frac{ C_{u}^{(d)}t^{\tf{4-d}{4}}}{\vepo^{d/4}}\ge\int_0^t\intrd\lab\KKSepsx\rab^2dx ds\ge\frac{ C_{l}^{(d)}t^{\tf{4-d}{4}}}{\vepo^{d/4}}
\eeq
\bpfs{Proof of \thmref{intense}} Let $T>0$ and $q\ge1$ be fixed and arbitrary, and let $(t,x)\in[0,T]\times\Rd$ and $d=1,2,3$.  
\ben\rencomrom 
\item
 Under the conditions \eqref{lcnd}, we have by \thmref{Holreg} a unique solution $U_{\vepo,\vept}$ to the L-KS SPDE \eqref{lksop} that is $L^{2q}(\Omega)$-bounded on $\T\times\Rd$, for every $q\ge1$.   Let $\mu^{t,x}_{\vepo}$ be the measure on $[0,t]\times\Rd$ defined by $$d\mu^{t,x}_{\vepo}(s,y)=\lab\KKSeptsxy\rab^{2}ds\,dy$$  and let
$|\mu^{t,x}_{\vepo}| = \mu^{t,x}_{\vepo}([0,t] \times \mathbb{R})$.  Taking the $2q$-th moment of the difference between our scaled L-KS SPDE and its deterministic counterpart---whose solution we denote by $u_{\vepo}$; using Burkholder's inequality followed by Jensen's inequality applied to the probability measure ${d\mu^{t,x}_{\vepo}(s,y)}/{|\mu^{t,x}_{\vepo}|}$; then using the linear growth condition ((b) in \eqref{lcnd}) on $a$, the $L^{2q}(\Omega)$-boundedness of $U_{\vepo,\vept}$ on $\T\times\Rd$, and the upper bound in \eqref{eplksl2},  we get
\begin{equation}
\begin{split}
&\mathbb{E}\left|U_{\vepo,\vept}(t,x)-u_{\vepo}(t,x)\right|^{2q}\\
&\le C\mathbb{E}\lab \int_{\Rd}\int_0^t \KKSeptsxy{\vept a(U_{\vepo,\vept}(s,y))} \sW(ds,dy)\rab^{2q}\\
&\le C\vept^{2q}\int_{\Rd}\int_0^t\E a^{2q}\lpa U_{\vepo,\vept}(s,y)\rpa
\frac{d\mu^{t,x}_{\vepo}(s,y)}{|\mu^{t,x}_{\vepo}|}|
\mu^{t,x}_{\vepo}|^{q}\\
&\le \frac{CT^{\tf{(4-d)}{4}q}\vept^{2q}}{\vepo^{(d/4)q}}\to0
\end{split}
\label{momineq1}
\end{equation}
as $\vepo,\vept$, and $\vept/\vepo^{d/8}$ approach $0$.
\item We prove it by contradiction.  So, assume there is a $T>0$ such that
\beq\lbl{contrass}
\lim_{\substack{\vepo,\vept\downarrow0\\
{\vept}/{{ \vepo}^{d/8}}\to\infty}}
\sup_{0\le s\le T}\sup_{{x\in\Rd}}\mathbb{E} U_{\vepo,\vept}^{2}(s,x)
<\infty;\ d=1,2,3
\eeq
and assume without loss of generality that $\un\equiv0$.  Observe that
\begin{equation}
\begin{split}
&\E\left|U_{\vepo,\vept}(t,x)\right|^2=\mathbb{E} \lab \int_{\Rd}\int_0^t \KKSeptsxy{\vept a(U_{\vepo,\vept}(s,y))} \sW(ds,dy)\rab^{2}\\
&= \vept^2\int_{\Rd}\int_0^t \lab\KKSeptsxy\rab^{2}
\E a^2\lpa U_{\vepo,\vept}(s,y)\rpa ds\,dy \\
&\ge K^2_{l}\vept^2\int_{\Rd}\int_0^t \lab\KKSeptsxy\rab^{2} ds\,dy
\ge \frac{\tilde{C}^{(d)}\vept^2t^{\tf{4-d}{4}}}{\vepo^{d/4}}; \quad d=1,2,3,
\end{split}
\label{as00}
\end{equation}
where we used the lower bound assumption $0<K_{l}\le a(u)$ and the lower bound in \eqref{eplksl2} to get the last two inequalities in \eqref{as00}.  Using the assumption in \eqref{contrass}, we arrive at the desired contradiction by taking the limit as $\vepo,\vept\searrow0$ in \eqref{as00} such that $\vept/\vepo^{d/8}\nearrow\infty$ for $d=1,2,3$.   
\een
The proof is complete.
\epfs
\section{Proof of \thmref{com}}
Here, we prove the change of measure transfer of properties from the canonical L-KS to nonlinear L-KS fourth order SPDEs, including the Swift-Hohenberg SPDE on subsets of $\{\Rp\times\Rd\}_{d=1}^{3}$.  

We say that a progressively measurable random field $X$ on the probability space
$\OFFtP$ satisfies Novikov's condition on $\T\times\S$ if
\begin{equation}\label{N}
\EP\left[\exp\left(\frac12\intop _{\T\times\S}X^2(t,x)dtdx\right)\right]<\infty.
\end{equation}
\bpfs{Proof of \thmref{com}}
\ben\rencomrom
\item (Transfer of law uniqueness)  We prove the more interesting direction (from zero to nonzero drift).  The proof of the reverse direction is similar and is omitted.  Suppose that uniqueness in law holds for the zero-drift L-KS SPDE $\eLKSazun$ and that
$\solViWiti,\ \filpspaceiti;\ i=1,2$ are solutions to the nonzero-drift L-KS SPDE $\eLKSabun$.  By assumption. we have 
\begin{equation}\label{L2com}
\Piti\Bigg[\int _{\T\times\S}{R_{V^{(i)}}^{2}}(t,x)dtdx<\infty \Bigg]=1,\ i=1,2.
\end{equation}
Define the sequence of stopping times $\{\tauni\}$ by
\begin{equation}\label{STimes}
\tauni\eqdef T\wedge\inf\left\{0\le t\le T;\int_{[0,t]\times\S}\RVi^2(s,x)ds dx=n\right\};
\ n\in\N,\ i=1,2,
\end{equation}
and let $\Wmi=\{\WitB,\Ft;0\le t\le T,B\in\BL\}$ be given by
$$\WitB\eqdef\WittB+\int _{[0,t]\times B}\RVi(s,x)dsdx;\ i=1,2.$$
Then, Novikov's condition \eqnref{N} and Girsanov's theorem for white noise (see
Theorem 2.2, Corollary 2.3, and Lemma 2.4 in \cite{Acom}) immediately gives us that $\Wmin=\{\WitBn,\Ft;0\le t\le T,B\in\BL\}$
is a white noise stopped at time $\tauni$, under the probability measure
$\Pni$ defined on $\FTi$ by 
\begin{equation*}
\frac{d\Pni}{d\Piti}=\RNbTtauniLYiWit;\ n\in\N,\ i=1,2,
\end{equation*}
where the Radon-Nikodym derivative is given by
\begin{equation*}
\begin{split}
&\RNbttauniBYiWit
\\ \eqdef&\exp\lbr\intop_{[0,t\wedge\tauni]\times B}-\lbk\RVi(s,x)\Wmit(ds,dx)
- \tf12\RVi^2(s,x) ds dx\rbk\rbr;
\end{split}
\end{equation*}
$0\le t\le T$, $B\in\BL$.  Consequently, $\solViWin$, $\filpspaceTni$ is a solution
to the zero-drift L-KS SPDE $\eLKSazun$ on $[0,T\wedge\tauni ]\times\S$ for each $i=1,2$ and $n\in\N$.
Clearly, for $i=1,2$,
\begin{equation}
\begin{split}
&\frac{d\Piti}{d\Pni}=\RNfTtauniLYiWi\\
\eqdef&\exp\lbr\intop_{[0,T\wedge\tauni]\times\S}\lbk\RVi(s,x)\Wmi(ds,dx)
- \tf12\RVi^2(s,x) ds dx\rbk\rbr;
\end{split}
\label{RDF}
\end{equation}
$n\in\N$.  Thus, for any set $\Lambda\in{\mathscr B}(\C(\T\times\S;\R))$
\begin{equation}
\begin{split}
\Poti\left[\Vo\in\Lambda,\tauno=T\right]
&=\EPno\left[1_{\{\Vo\in\Lambda,\tauno=T\}}{\RNfTtaunoLYoWo}\right]\\
&=\EPnt\left[1_{\{\Vt\in\Lambda,\taunt=T\}}{\RNfTtauntLYtWt}\right]\\
&=\Ptti\left[\Vt\in\Lambda,\taunt=T\right];\ \forall n\in\N,
\end{split}
\label{uil}
\end{equation}
where we have used the uniqueness in law assumption on $\eLKSazun$ (comparing
the $\Vi$'s only on $\Omega_n^{(i)}\eqdef\{\tauni=T\}$ for each $n$), \eqnref{STimes},
and \eqnref{RDF} to get the second equality in \eqnref{uil}.
By \eqnref{L2com} and \eqnref{STimes}, we obtain that
$\lim \nolimits _{n\rightarrow \infty }\Piti[\tauni=T]=1\mbox{ for } i=1,2$.
Thus, taking the limit as $n\rightarrow \infty $ in \eqnref{uil} yields that the
law of $\Vo$ under $\Poti$ is the same as that of $\Vt$ under $\Ptti$.  I.e.,
we have uniqueness in law for the non-zero drift L-KS SPDE $\eLKSabun$.  
\item  (Law equivalence)  Let $(V,\sW^{(1)})$ be a solution $($weak or strong\/$)$ to the nonzero-drift L-KS SPDE $\eLKSabun$ on $\filqspace$; and let $(U,\sW^{(2)})$ be a solution $($weak or strong\/$)$ to the zero-drift L-KS SPDE $\eLKSazun$ on $\filpspace$.  Then, uniqueness in law for the L-KS SPDE $\eLKSabun$ follows from the uniqueness in law assumption for $\eLKSazun$, the almost sure $L^2(\T\times\S,\lambda)$ condition on $\RV$, and part (i) of \thmref{com}.

Replacing $\Vi$ in \eqref{STimes} by $U$ and then $V$, we get the definitions of the stopping times sequences $\{\taunU\}$ and $\{\taunV\}$, respectively.  Let $\Wtim=\{\WtitB,\Ft;0\le t\le T,B\in\BR\}$
be given by
$$\WtitB\eqdef\WttB-\int _{[0,t]\times B}\RU(s,x)dsdx.$$
Then, Novikov's condition and Girsanov's theorem for white noise (see
Theorem 2.2, Corollary 2.3, and Lemma 2.4 in \cite{Acom})
immediately give us that, for $n\in\N$, $\Wtmn=\{\WtitBn,\Ft;0\le t\le T,B\in\BR\}$ is a white noise stopped at time $\taunU$, under the
probability measure $\Ptn$ defined on $\FT$ by the Radon-Nikodym derivative
\begin{equation}\label{RNf}
\begin{split}
&\frac{d\Ptn}{d\P}=\RNfTtaunWt
\\&\eqdef\ \exp\lbr\intop_{[0,T\wedge\taunU]\times\S}\RU(s,x)\right.\Wm^{(2)}(ds,dx)
\left.- \frac12\intop_{[0,T\wedge\taunU]\times\S} \RU^{2}(s,x) ds dx\rbr.
\end{split}
\end{equation}
Thus, $\solUWt$, $\filTpspace$ is a solution to the nonzero-drift L-KS SPDE $\eLKSabun$ on $[0,T\wedge\taunU]\times\S$, for each $n\in\N$.  As a result, for any set $\Lambda\in{\sB\lpa\C(\T\times\S;\R)\rpa}$ we get
\begin{equation}
\begin{split}
\Q[V\in\Lambda,\taunV=T]&=\Ptn[U\in\Lambda,\taunU=T]
\\&=\EP\left[1_{\{U\in\Lambda,\taunU=T\}}\RNfTtaunW\right];\quad n\in\N.
\end{split}
\label{1}
\end{equation}
To see \eqref{1} note that, on the event $\Omega_n^U\eqdef\{\omega\in\Omega^{(2)};\taunU(\omega)=T\}$,
$\solUWt$ is a solution to $\eLKSabun$ on $\T\times\S$, under $\Ptn$.  Thus, the first equality in \eqref{1} follows from the uniqueness in law for $\eLKSabun$ and the definitions of $\taunU$ and $\taunV$.
By the $L^2$ assumption on $\RV$ and the definition of $\taunV$, we have $\lim_{n\to\infty} \Q[\taunV=T]=1$; so, taking limits in \eqref{1} we get
\begin{equation}
\Q[V\in\Lambda]=\lim_{n\to\infty}\Ptn[U\in\Lambda,\taunU=T]
=\lim_{n\to\infty}\EP\left[1_{\{U\in\Lambda,\taunU=T\}}\RNfTtaunWt\right].
\label{2}
\end{equation}
Obviously, if $\P[U(\cdot,\cdot)\in\Lambda]=0$ then $\EP\left[1_{\{U(\cdot,\cdot)\in\Lambda,\taunU=T\}}\RNfTtaunWt\right]=0$ for each $n$; thus,
$$\Q[V(\cdot,\cdot)\in\Lambda]=\lim_{n\to\infty}\EP\left[1_{\{U(\cdot,\cdot)\in\Lambda,\taunU=T\}}\RNfTtaunWt\right]=0.$$
I.e., $\LawVQ$ is absolutely continuous with respect to $\LawUP$ on $\sB\lpa\C(\T\times\S;\R)\rpa$.   The absolute continuity of $\LawUP$ with respect to $\LawVQ$ is proved by a similar argument, and we omit it.
\een
The proof of the last part of \thmref{com} follows since the square of the drift/diffusion ratios given by the random fields
\beq\lbl{ratcondver}
\bsp
R^{2}_{U}(t,x)=\sum_{k=0}^{2N}\tilde{c}_{k}U^{k}(t,x)\mbox{ and }R^{2}_{V}(t,x)=\sum_{k=0}^{2N}\tilde{c}_{k}V^{k}(t,x)\mbox{ for }\tilde{c}_{i}\in\R  
\end{split}
\eeq
are continuous and thus almost surely bounded on the compact set $[0,T]\times\S$.  Thus, $\RU$ and $\RV$ are in $L^2(\T\times\S,\lambda)$, almost surely, and parts (i) and (ii) of \thmref{com} follow in the case $a(u)=\kappa\in\R\setminus\{0\}$,  $b(u)=\sum_{k=0}^{N}c_{k}u^{k}$ for $c_{k}\in\R$, $k=0,\ldots,N$, and $N\ge0$, and $u_0\in \mathrm{C}_{c}^{2,\gamma}(\S;\R)\mbox{ and nonrandom}$..  The uniqueness assertion for $\eLKSabun$ follows from (i) together with \lemref{pathuniqonS} and the H\"older equivalence assertion follows from (ii).
\epfs
\appendix
\section{Uniqueness lemma}\lbl{A}
Throughout this Appendix, we reserve the notation $\S$ solely for the $d$-dimensional rectangle $\prod_{i=1}^{d}[0,L_{i}]$, $d=1,2,3$; and, as before, $\T=[0,T]$ for some fixed but arbitrary $T>0$.  We now prove pathwise uniqueness for $\eLKSazun$ with Dirichlet boundary conditions.  The following lemma is useful for the last part of \thmref{com}.
\blm\lbl{pathuniqonS}
Pathwise uniqueness, and hence uniqueness in law, holds for the zero-drift L-KS SPDE $\eLKSazun$ on $\Rp\times\S$ whenever $a\equiv \kappa\neq0$.
\elm
\bpf  Assume without loss of generality that $d=1$, $\S=[0,1]$, and fix $\vep>0,\vth\in\R$.  Assume further that $(U^{(1)},\sW)$ and $(U^{(2)},\sW)$ are two continuous solutions to $\eLKSazun$ on $\Rp\times\S$ on the same usual probability space $\OFFtP$ and with respect to the same space-time white noise $\sW$.  Let $D(t,x)=U^{(1)}(t,x)-U^{(2)}(t,x)$ for $(t,x)\in\Rp\times\S$,  
and let $$\Phi^{\infty}_{\mbox\tiny{Dir}}(\S;\R):=\lbr\vph\in\C^{\infty}(\Rd;\R); \varphi=\D\varphi=0\mbox{ on }\pa\S\rbr.$$ 

Now, the continuous difference random field $D$ satisfies 
\begin{equation}\label{integral1}
\bsp
&\int_0^1D(t,x)\vph(x) dx
\\&=-\frac\vep8\int_0^t\int_0^1D(s,x) \lpa\vph^{(4)}(x)+4\vth\vph^{(2)}(x)+4\vth^{2}\vph(x)\rpa dx ds
\end{split}
\eeq
for every $\varphi\in\Phi^{\infty}_{\mbox\tiny{Dir}}(\S;\R)$, $t\in\T,\ \mbox{a.s. } \P.$
Manifestly, this implies that $D(t,x)=0$ on $[0,T]\times[0,1]$ a.s. $\P$.
To see this, choose $\varphi_m(x)=\sin(m\pi x)$ for $m\in\Bbb{N}$ and let 
$$C_m(t):=\int_0^1D(t,x)\varphi_m(x)dx;\ \forall m\in\Bbb{N}, t\in\T.$$ 
So, by \eqref{integral1}, we have
\begin{equation}
C_m(t)=\lpa-\frac{\vep m^4\pi^4}{8}+\frac{\vep\vth m^2\pi^2}{2}-\frac{\vep\vth^{2}}{2}\rpa\int_0^tC_m(s)ds; 
\ \mbox{a.s. } \P,
\label{integral2}
\end{equation}
for all $t\in[0,T]$ and $m\in\Bbb{N}$, which obviously implies that, for each $m$, $C_m(t)=0$ for all $t\in\T$ a.s. $\P$.   
Now, since all the Fourier Sine coefficients, $C_m(t)$, for $D(t,x)$ are zero and 
the continuous solutions of the SPDE $\eLKSazun$ vanish at $x=0$ and $x=1$ for all $t$
(and hence $D(t,0)=D(t,1)=0$ $\forall t$), we get that $D(t,x)=0$ 
on $[0,T]\times[0,1]$ a.s. $\P$.  The arbitrariness of $T$ now completes the proof.
\epf
\section{Frequent acronyms and notations key}\lbl{B}
\begin{enumerate}\renewcommand{\labelenumi}{\Roman{enumi}.}
\item {\textbf{Acronyms}}\vspace{2mm}
\begin{enumerate}\renewcommand{\labelenumii}{(\arabic{enumii})}
\item BTBM: Brownian-time Brownian motion.
\item BTBM SIE: BTBM stochastic integral equation.
\item IBTBAP: imaginary Brownian-time Brownian angle process.
\item KS: Kuramoto Sivashinsky.
\item SH: Swift-Hohenberg.
\end{enumerate}
\vspace{2.5mm}
\item {\textbf{Notations}}\vspace{2mm}
\begin{enumerate}\renewcommand{\labelenumii}{(\arabic{enumii})}
\item $\KKSepthtx$ the $(\vep,\vth)$ L-KS kernel (see \eqref{vepvthLKS}).
\item $\KKStx=\KKSepthootx$ the canonical (or $(1,1)$) L-KS kernel (see \eqref{ibtbapkernel} and \eqref{vepvthLKS}).
\item $\Ksfotx:=\KKSepthoztx$ the zero-angle canonical (or simple fourth order) kernel (see \eqref{sfoker} and \eqref{vepvthLKS}).
\item $\C^{k,\gamma}$ the set of functions with $\gamma$-H\"older continuous $k$-th derivative, $k=0,1,2,\ldots$ 
\item $\C^{k,\gamma}_{b}$ the set of functions with $\gamma$-H\"older continuous and bounded $k$-th derivative, $k=0,1,2,\ldots$ 
\item $\T:=[0,T]$ for some fixed and arbitrary $T>0$
\item $M_{p}(t):=\sup_{x\in\Rd}\E\lab U(t,x)\rab^{p}$
\end{enumerate}
\end{enumerate}

\end{document}